\documentclass[12pt,reqno]{amsart}
\usepackage{mathrsfs}
\usepackage{cases}
\usepackage{epic}
\usepackage{amsfonts}
\usepackage{tikz}
\usepackage{tikz-cd}
\usepackage{graphicx}
\usepackage{amsmath}
\usepackage{amssymb, upgreek}
\usepackage{bm}
\usepackage{multirow,longtable,dcolumn} 
\usepackage{latexsym,todonotes}
\usepackage{pdflscape}
\usepackage[all]{xypic}
\usepackage[all]{xy}
\usepackage{color}
\usepackage{colordvi}
\usepackage{multicol}
\usepackage{mathtools}
\usepackage[normalem]{ulem}
\usepackage[linktocpage=true]{hyperref}
\hypersetup{colorlinks,linkcolor=blue,urlcolor=cyan,citecolor=blue}
\usepackage{graphicx}
\NewDocumentCommand{\sslash}{s}{%
	\IfBooleanTF{#1}
	{\big/\mkern-7mu\big/}
	{/\mkern-6mu/}%
}
\makeatletter
\newsavebox{\@brx}
\newcommand{\llangle}[1][]{\savebox{\@brx}{\(\m@th{#1\langle}\)}%
	\mathopen{\copy\@brx\kern-0.5\wd\@brx\usebox{\@brx}}}
\newcommand{\rrangle}[1][]{\savebox{\@brx}{\(\m@th{#1\rangle}\)}%
	\mathclose{\copy\@brx\kern-0.5\wd\@brx\usebox{\@brx}}}
\makeatother

\tikzcdset{scale cd/.style={every label/.append style={scale=#1},
		cells={nodes={scale=#1}}}}
\DeclareMathAlphabet{\mathbbb}{U}{bbold}{m}{n}
\DeclareMathOperator{\dimv}{\underline{\dim}}

\makeatletter

\newcommand{\Rmnum}[1]{\textup{\expandafter\@slowromancap\romannumeral #1@}}
\makeatother

\allowdisplaybreaks
\def \wt{\mathrm{wt}}
\def \tUU{(\tU\otimes\tU)}
\def \tUUi {(\tU\otimes\tU)^\imath}

\begin{document}
	\input xy
	\xyoption{all}
	\newcommand{\iLa}{\Lambda^{\imath}}
	\newcommand{\iadd}{\operatorname{iadd}\nolimits}
	\renewcommand{\mod}{\operatorname{mod}\nolimits}
	\newcommand{\fproj}{\operatorname{f.proj}\nolimits}
	\newcommand{\Fac}{\operatorname{Fac}\nolimits}
	\newcommand{\ci}{{\I}_{\btau}}
	\newcommand{\proj}{\operatorname{proj}\nolimits}
	\newcommand{\inj}{\operatorname{inj}\nolimits}
	\newcommand{\rad}{\operatorname{rad}\nolimits}
	\newcommand{\Span}{\operatorname{Span}\nolimits}
	\newcommand{\soc}{\operatorname{soc}\nolimits}
	\newcommand{\ind}{\operatorname{inj.dim}\nolimits}
	\newcommand{\Ginj}{\operatorname{Ginj}\nolimits}
	\newcommand{\res}{\operatorname{res}\nolimits}
	\newcommand{\np}{\operatorname{np}\nolimits}
	\newcommand{\Mor}{\operatorname{Mor}\nolimits}
	\newcommand{\Mod}{\operatorname{Mod}\nolimits}
	\newcommand{\End}{\operatorname{End}\nolimits}
	\newcommand{\lf}{\operatorname{l.f.}\nolimits}
	\newcommand{\Iso}{\operatorname{Iso}\nolimits}
	\newcommand{\Aut}{\operatorname{Aut}\nolimits}
	\newcommand{\Rep}{\operatorname{Rep}\nolimits}
	
	\newcommand{\colim}{\operatorname{colim}\nolimits}
	\newcommand{\gldim}{\operatorname{gl.dim}\nolimits}
	\newcommand{\cone}{\operatorname{cone}\nolimits}
	\newcommand{\rep}{\operatorname{rep}\nolimits}
	\newcommand{\Ext}{\operatorname{Ext}\nolimits}
	\newcommand{\Tor}{\operatorname{Tor}\nolimits}
	\newcommand{\Hom}{\operatorname{Hom}\nolimits}
	\newcommand{\Top}{\operatorname{top}\nolimits}
	\newcommand{\Coker}{\operatorname{Coker}\nolimits}
	\newcommand{\thick}{\operatorname{thick}\nolimits}
	\newcommand{\rank}{\operatorname{rank}\nolimits}
	\newcommand{\Gproj}{\operatorname{Gproj}\nolimits}
	\newcommand{\Len}{\operatorname{Length}\nolimits}
	\newcommand{\RHom}{\operatorname{RHom}\nolimits}
	\renewcommand{\deg}{\operatorname{deg}\nolimits}
	\renewcommand{\Im}{\operatorname{Im}\nolimits}
	\newcommand{\Ker}{\operatorname{Ker}\nolimits}
	\newcommand{\Coh}{\operatorname{Coh}\nolimits}
	\newcommand{\Id}{\operatorname{Id}\nolimits}
	\newcommand{\Qcoh}{\operatorname{Qch}\nolimits}
	\newcommand{\CM}{\operatorname{CM}\nolimits}
	\newcommand{\sgn}{\operatorname{sgn}\nolimits}
	\newcommand{\utMH}{\operatorname{\cm\ch(\iLa)}\nolimits}
	\newcommand{\GL}{\operatorname{GL}}
	\newcommand{\Perv}{\operatorname{Perv}}
	
	\newcommand{\IC}{\operatorname{IC}}
	\def \hU{\widehat{\U}}
	\def \hUi{\widehat{\U}^\imath}
	\newcommand{\bb}{\psi_*}
	\newcommand{\bvs}{{\boldsymbol{\varsigma}}}
	\def \ba{\mathbf{a}}
	\newcommand{\vs}{\varsigma}
	\def \bfk {\mathbf{k}}

	\def \bd{\mathbf{d}}
	\newcommand{\e}{{\bf 1}}
	\newcommand{\EE}{E^*}
	\newcommand{\dbl}{\operatorname{dbl}\nolimits}
	\newcommand{\ga}{\gamma}
	\newcommand{\tM}{\cm\widetilde{\ch}}
	\newcommand{\la}{\lambda}
	
	\newcommand{\For}{\operatorname{{\bf F}or}\nolimits}
	\newcommand{\coker}{\operatorname{Coker}\nolimits}
	\newcommand{\rankv}{\operatorname{\underline{rank}}\nolimits}
	\newcommand{\diag}{{\operatorname{diag}\nolimits}}
	\newcommand{\swa}{{\operatorname{swap}\nolimits}}
	\newcommand{\supp}{{\operatorname{supp}}}
	
	\renewcommand{\Vec}{{\operatorname{Vec}\nolimits}}
	\newcommand{\pd}{\operatorname{proj.dim}\nolimits}
	\newcommand{\gr}{\operatorname{gr}\nolimits}
	\newcommand{\id}{\operatorname{id}\nolimits}
	\newcommand{\aut}{\operatorname{Aut}\nolimits}
	\newcommand{\Gr}{\operatorname{Gr}\nolimits}
	
	\newcommand{\pdim}{\operatorname{proj.dim}\nolimits}
	\newcommand{\idim}{\operatorname{inj.dim}\nolimits}
	\newcommand{\Gd}{\operatorname{G.dim}\nolimits}
	\newcommand{\Ind}{\operatorname{Ind}\nolimits}
	\newcommand{\add}{\operatorname{add}\nolimits}
	\newcommand{\pr}{\operatorname{pr}\nolimits}
	\newcommand{\oR}{\operatorname{R}\nolimits}
	\newcommand{\oL}{\operatorname{L}\nolimits}
	\def \brW{\mathrm{Br}(W_\btau)}
	\newcommand{\Perf}{{\mathfrak Perf}}
	\newcommand{\cc}{{\mathcal C}}
	\newcommand{\gc}{{\mathcal GC}}
	\newcommand{\ce}{{\mathcal E}}
	\newcommand{\calI}{{\mathcal I}}
	\newcommand{\cs}{{\mathcal S}}
	\newcommand{\cf}{{\mathcal F}}
	\newcommand{\cx}{{\mathcal X}}
	\newcommand{\cy}{{\mathcal Y}}
	\newcommand{\ct}{{\mathcal T}}
	\newcommand{\cu}{{\mathcal U}}
	\newcommand{\cv}{{\mathcal V}}
	\newcommand{\cn}{{\mathcal N}}
	\newcommand{\mcr}{{\mathcal R}}
	\newcommand{\ch}{{\mathcal H}}
	\newcommand{\ca}{{\mathcal A}}
	\newcommand{\cb}{{\mathcal B}}
	\newcommand{\cj}{{\mathcal J}}
	\newcommand{\cl}{{\mathcal L}}
	\newcommand{\cm}{{\mathcal M}}
	\newcommand{\cp}{{\mathcal P}}
	\newcommand{\cg}{{\mathcal G}}
	\newcommand{\cw}{{\mathcal W}}
	\newcommand{\co}{{\mathcal O}}
	\newcommand{\cq}{{\mathcal Q}}
	\newcommand{\cd}{{\mathcal D}}
	\newcommand{\ck}{{\mathcal K}}
	\newcommand{\calr}{{\mathcal R}}
	\newcommand{\cz}{{\mathcal Z}}
	\newcommand{\ol}{\overline}
	\newcommand{\ul}{\underline}
	\newcommand{\st}{[1]}
	\newcommand{\ow}{\widetilde}
	\renewcommand{\P}{\mathbf{P}}
	\newcommand{\pic}{\operatorname{Pic}\nolimits}
	\newcommand{\Spec}{\operatorname{Spec}\nolimits}
	\newcommand{\Fr}{\mathrm{Fr}}
	\newcommand{\Gp}{\mathrm{Gp}}
	\newcommand{\fU}{\mathfrak{U}}
	
	\newtheorem{innercustomthm}{{\bf Theorem}}
	\newenvironment{customthm}[1]
	{\renewcommand\theinnercustomthm{#1}\innercustomthm}
	{\endinnercustomthm}
	
	\newtheorem{innercustomcor}{{\bf Corollary}}
	\newenvironment{customcor}[1]
	{\renewcommand\theinnercustomcor{#1}\innercustomcor}
	{\endinnercustomthm}
	
	\newtheorem{innercustomprop}{{\bf Proposition}}
	\newenvironment{customprop}[1]
	{\renewcommand\theinnercustomprop{#1}\innercustomprop}
	{\endinnercustomthm}
	
	\newtheorem{theorem}{Theorem}[section]
	\newtheorem{acknowledgement}[theorem]{Acknowledgement}
	\newtheorem{algorithm}[theorem]{Algorithm}
	\newtheorem{axiom}[theorem]{Axiom}
	\newtheorem{case}[theorem]{Case}
	\newtheorem{claim}[theorem]{Claim}
	\newtheorem{conclusion}[theorem]{Conclusion}
	\newtheorem{condition}[theorem]{Condition}
	\newtheorem{conjecture}[theorem]{Conjecture}
	\newtheorem{construction}[theorem]{Construction}
	\newtheorem{corollary}[theorem]{Corollary}
	\newtheorem{criterion}[theorem]{Criterion}
	\newtheorem{definition}[theorem]{Definition}
	\newtheorem{example}[theorem]{Example}
	\newtheorem{assumption}[theorem]{Assumption}
	\newtheorem{lemma}[theorem]{Lemma}
	\newtheorem{notation}[theorem]{Notation}
	\newtheorem{problem}[theorem]{Problem}
	\newtheorem{proposition}[theorem]{Proposition}
	\newtheorem{solution}[theorem]{Solution}
	\newtheorem{summary}[theorem]{Summary}
	\newtheorem{hypothesis}[theorem]{Hypothesis}
	\newtheorem*{thm}{Theorem}
	
	\theoremstyle{remark}
	\newtheorem{remark}[theorem]{Remark}
	
	\def \Br{\mathrm{Br}}
	\newcommand{\tK}{K}
	
	\newcommand{\tk}{\widetilde{k}}
	\newcommand{\tU}{\widetilde{{\mathbf U}}}
	\newcommand{\Ui}{{\mathbf U}^\imath}
	\newcommand{\tUi}{\widetilde{{\mathbf U}}^\imath}
	\newcommand{\qbinom}[2]{\begin{bmatrix} #1\\#2 \end{bmatrix} }
	\newcommand{\ov}{\overline}
	\newcommand{\tMHg}{\operatorname{\widetilde{\ch}(Q,\btau)}\nolimits}
	\newcommand{\tMHgop}{\operatorname{\widetilde{\ch}(Q^{op},\btau)}\nolimits}
	
	\newcommand{\rMHg}{\operatorname{\ch_{\rm{red}}(Q,\btau)}\nolimits}
	\newcommand{\dg}{\operatorname{dg}\nolimits}
	\def \fu{{\mathfrak{u}}}
	\def \fv{{\mathfrak{v}}}
	\def \sqq{{\mathbbb{v}}}
	\def \bp{{\mathbf p}}
	\def \bv{{\mathbf v}}
	\def \bw{{\mathbf w}}
	\def \bA{{\mathbf A}}
	\def \bL{{\mathbf L}}
	\def \bF{{\mathbf F}}
	\def \bS{{\mathbf S}}
	\def \bC{{\mathbf C}}
	\def \bU{{\mathbf U}}
	\def \U{{\mathbf U}}
	\def \btau{\varpi}
	\def \La{\Lambda}
	\def \Res{\Delta}
	\newcommand{\ev}{\bar{0}}
	\newcommand{\odd}{\bar{1}}
	\def \fk{\mathfrak{k}}
	\def \ff{\mathfrak{f}}
	\def \fp{{\mathfrak{P}}}
	\def \fg{\mathfrak{g}}
	\def \fn{\mathfrak{n}}
	\def \gr{\mathfrak{gr}}
	\def \Z{\mathbb{Z}}
	\def \F{\mathbb{F}}
	\def \D{\mathbb{D}}
	\def \C{\mathbb{C}}
	\def \N{\mathbb{N}}
	\def \Q{\mathbb{Q}}
	\def \G{\mathbb{G}}
	\def \P{\mathbb{P}}
	\def \K{\mathbb{K}}
	\def \E{\mathbb{K}}
	\def \I{\mathbb{I}}
	
	\def \eps{\varepsilon}
	\def \BH{\mathbb{H}}
	\def \btau{\varrho}
	\def \cv{\varpi}
	
	\def \tR{\widetilde{\bf R}}
	\def \tRZ{\widetilde{\bf R}_\cz}
	\def \hR{\widehat{\bf R}}
	\def \hRZ{\widehat{\bf R}_\cz}
	\def\tRi{\widetilde{\bf R}^\imath}
	\def\hRi{\widehat{\bf R}^\imath}
	\def\tRiZ{\widetilde{\bf R}^\imath_\cz}
	\def\reg{\mathrm{reg}}
	
	\def\hRiZ{\widehat{\bf R}^\imath_\cz}
	\def \tTT{\widetilde{\mathbf{T}}}
	\def \TT{\mathbf{T}}
	\def \br{\mathbf{r}}
	\def \bp{{\mathbf p}}
	\def \tS{\texttt{S}}
	\def \bq{{\bm q}}
	\def \bvt{{v}}
	\def \bs{{ r}}
	\def \tt{{v}}
	\def \k{k}
	\def \bnu{\bm{\nu}}
	\def\bc{\mathbf{c}}
	\def \ts{\textup{\texttt{s}}}
	\def \tt{\textup{\texttt{t}}}
	\def \tr{\textup{\texttt{r}}}
	\def \tc{\textup{\texttt{c}}}
	\def \tg{\textup{\texttt{g}}}
	\def \bW{\mathbf{W}}
	\def \bV{\mathbf{V}}

	\newcommand{\browntext}[1]{\textcolor{brown}{#1}}
	\newcommand{\greentext}[1]{\textcolor{green}{#1}}
	\newcommand{\redtext}[1]{\textcolor{red}{#1}}
	\newcommand{\bluetext}[1]{\textcolor{blue}{#1}}
	\newcommand{\brown}[1]{\browntext{ #1}}
	\newcommand{\green}[1]{\greentext{ #1}}
	\newcommand{\red}[1]{\redtext{ #1}}
	\newcommand{\blue}[1]{\bluetext{ #1}}
	\numberwithin{equation}{section}
	\renewcommand{\theequation}{\thesection.\arabic{equation}}
	
	\newcommand{\wtodo}{\rightarrowdo[inline,color=orange!20, caption={}]}
	\newcommand{\lutodo}{\rightarrowdo[inline,color=green!20, caption={}]}
	\def \tT{\widetilde{\mathcal T}}
	
	\def \tTL{\tT(\iLa)}
	\def \iH{\widetilde{\ch}}
	\newcommand{\EuFo}[1]{\left \langle #1 \right \rangle}
	
	\title[Quantum symmetric pairs via Hall algebras]{Quantum symmetric pairs via Hall algebras}

	\author[Ming Lu]{Ming Lu}
	\address{Department of Mathematics, Sichuan University, Chengdu 610064, P.R.China}
	\email{luming@scu.edu.cn}

	\author[Zhuoyi Zhao]{Zhuoyi Zhao}
	\address{Department of Mathematics, Nanjing University, Nanjing 210093, P.R.China}
	\email{502023210029@smail.nju.edu.cn}

	
	\maketitle
	\begin{quote}\begin{center}
			{\em Dedicated to Professor Bin Zhu on the occasion of his 60th birthday}
		\end{center}
	\end{quote}
	
	\begin{abstract}
		A quantum symmetric pair consists of a quantum group
		$\tU$ and its coideal subalgebra $\tUi$. The Hall algebra constructions of $\tU$ and
		$\tUi$ are given by Bridgeland and Lu--Wang, respectively. In this paper, we construct a Hall algebra framework for the coideal subalgebra structure of $\tUi$ in $\tU$, and for the quantum symmetric pair $(\tU,\tUi)$. As an application, we prove that the natural embedding $\imath:\tUi\to \tU$, and the coproduct $\Delta:\tUi\to \tUi\otimes \tU$ preserve the integral forms of $\tUi$ and $\tU$,  which are used to construct the dual canonical bases.
	\end{abstract}

	\setcounter{tocdepth}{1}
	\tableofcontents

	\section{Introduction}
	
	Let $\U=\U_v(\fg)$ be the Drinfeld--Jimbo quantum group; see \cite{Dr87,J85}. Ringel \cite{Rin90} introduced the Hall algebra on quiver representations to realize the positive part $\U^+$ of $\U$ of Dynkin type; see \cite{Gr95} for an extension to Kac-Moody type. Inspired by Ringel's construction, Lusztig gave a geometric realization of $\U^+$ by using perverse sheaves over quiver varieties, who also constructed the canonical basis on $\U^+$ \cite{Lus90,Lus93}.  
	
	It is a natural question to realize the whole quantum group $\U$ by using Hall algebras. Many efforts are taken to solve this question, see \cite{X97,Ka98,PX00} for some earlier attempts, and see \cite{T06,XX08} for the construction of derived Hall algebras; Bridgeland \cite{Br13} solved this question by considering a certain Hall algebra of the $2$-periodic complexes of projective modules, in fact, he realized the Drinfeld double $\tU$, a variant of $\U$ with the Cartan part doubled. Bridgeland's construction has found further generalization and improvements; see \cite{Gor18} for semi-derived Hall algebras of Frobenius categories, and see \cite{LP21} for semi-derived Ringel-Hall algebras of hereditary categories.

	Inspired by Bridgeland's Hall algebras, Qin \cite{Qin} realized the whole quantum group $\tU$ of type ADE by perverse sheaves over some cyclic quiver variety, based on the construction in \cite{Lus90,Na01,HL15}. Qin's geometric realization gives an integral and positive basis on $\tU$, which contains as subsets
	the (mildly rescaled) dual canonical bases of Lusztig for $\U^+$ and $\U^-$.

	The iquantum group $\Ui$ is introduced in \cite{Let99} as a coideal subalgebra of the quantum group $\U$, with the Satake diagram as input; see \cite{Ko14} for general iquantum groups of Kac-Moody type. Following \cite{Let99,Ko14}, we call $(\tU,\tUi)$ and $(\U,\Ui)$ quantum symmetric pairs. iQuantum groups are a vast generalization of quantum groups, since quantum groups can be viewed as iquantum groups of diagonal type. A breakthrough on iquantum groups is the discovery of icanonical bases by Bao and Wang \cite{BW18,BW18b}.
	In this paper, we focus on quasi-split iquantum groups of type ADE.  
	
	A Dynkin iquiver $(Q,\varrho)$ consists of a Dynkin quiver $Q$ and an involution $\varrho$ of $Q$. Associated to $(Q,\varrho)$, the iquiver algebra $\Lambda^\imath$ is defined in \cite{LW22a}. The (twisted) semi-derived Ringel-Hall algebra (also called iHall algebra) $\widetilde{\ch}(\bfk Q,\varrho)$ of $\Lambda^\imath$ provides a categorical realization of the universal iquantum group $\tUi$. Compared with $\Ui$, the universal version $\tUi$ has various central elements, and $\Ui$ can be reproduced from $\tUi$ by a central reduction.  The quantum group $\tU$ is also reproduced by using the iquiver of  diagonal type. In the geometric setting, based on \cite{Qin,Sch19}, we \cite{LW21b} use perverse sheaves on NKS categories/varieties of $\Lambda^\imath$ to realize $\tUi$ of type ADE. Consequently, an integral and positive basis is also constructed on $\tUi$, called the dual canonical basis.  
	
	Recently, the first author and Pan \cite{LP25} gave a new construction of the dual canonical basis of $\tUi$ of type ADE via rescaled iHall basis and Lusztig's lemma, connecting earlier constructions in \cite{LW22a} and \cite{LW21b}. 
	In order to construct  and study the dual canonical bases for non-simpy-laced Dynkin type, the framework of iHopf algebras is introduced in \cite{CLPRW,CLPRW2}, as a generalization of the Drinfeld double construction.
	
	The goal of this paper is to realize the natural embedding $\imath:\tUi\to \tU$ and the coideal subalgebra structure of $\tUi$ in $\tU$ in the Hall algebra setting, solving the open question raised in \cite[\S11.1]{LW24}. The strategy is to compare and unify the realizations of $\tU$ and $\tUi$ using iHall algebras and iHopf algebras. As an application, this enables us to compare the integral forms and dual canonical bases of $\tU$ and $\tUi$.
	Let us summarize the main results in the following.
	
	For a Dynkin iquiver $(Q,\varrho)$, as proved in \cite{LW22a}, there exists an algebra isomorphism $\widetilde{\psi}:\tUi\rightarrow \widetilde{\ch}(Q,\varrho)$, where $\widetilde{\ch}(Q,\varrho)$ is the generic iHall algebra; cf. Lemma \ref{lem:Hall-iQG}. For the iquiver of diagonal type $(Q^{\dbl},\swa)$ (see Example \ref{ex:diagquiver}), there exists an algebra isomorphism $\widetilde{\psi}^{\dbl}:\tU\rightarrow \widetilde{\ch}(Q^{\dbl},\swa)$; see \cite{Br13}.

	The classic extended Hall algebra defined in \cite{Rin90,Gr95} can be viewed as a subalgebra of $\widetilde{\ch}(Q^{\dbl},\swa)$, which is denoted by $\widetilde{\mathbf{B}}(Q)$; see \S\ref{subsec:iHall-iHopf}. Then $\widetilde{\mathbf{B}}(Q)$ is a Hopf algebra with a Hopf pairing $\varphi(-,-)$; see \cite{Gr95}. This enables us to apply the machinery of iHopf algebras introduced in \cite{CLPRW} to  formulate the iHopf algebra $\widetilde{\mathbf{B}}^\imath_\varrho(Q)$, and the iHopf algebra of diagonal type $D(\widetilde{\mathbf{B}}(Q))$. 
	It is proved in Propositions \ref{prop:iHall-ihopf}, \ref{prop: iso of diagonal type} that there exist algebra isomorphisms
	\begin{align*}
		\widetilde{\Phi}: \widetilde{\ch}(Q,\varrho)\longrightarrow \widetilde{\mathbf{B}}^\imath_\varrho(Q),\qquad \widetilde{\Phi}^{\dbl}: \widetilde{\ch}(Q^{\dbl},\swa)\longrightarrow D(\widetilde{\mathbf{B}}( Q)).        \end{align*}
	We transfer the Hopf algebra structure of $D(\widetilde{\mathbf{B}}(Q))$ given in \cite[Proposition 2.7]{CLPRW} to $\widetilde{\ch}(Q^{\dbl},\swa)$, and prove that $\widetilde{\psi}^{\dbl}:\tU\rightarrow \widetilde{\ch}(Q^{\dbl},\swa)$ is an isomorphism as Hopf algebras, strengthening Bridgeland's result \cite{Br13}; see Theorem \ref{thm:Hopf-iHall} and Corollary \ref{cor:Bridgeland-Hopf}.

	We prove in Theorem \ref{thm:embedding} that there exists an injective algebra homomorphism $\widetilde{\Omega}: \widetilde{\ch}(Q,\varrho)\longrightarrow \widetilde{\ch}(Q^{\dbl},\swa)$ such that
	$$ \xymatrix{ \tUi \ar[rr]^{\imath }\ar[d]^{\widetilde{\psi}} && \tU \ar[d]^{\widetilde{\psi}^{\dbl}} 
		\\
		\widetilde{\ch}(Q,\varrho) \ar[rr]^{\widetilde{\Omega}} && \widetilde{\ch}(Q^{\dbl},\swa)} $$
	commutes. This result gives a realization of the quantum symmetric pair $(\tU,\tUi)$ in the Hall algebra setting. The technique is based on \cite[Proposition 2.12]{CLPRW}, and the crucial point/harder part is to find the closed formula of the   $\varrho$-twisted compatible linear map $\chi:\widetilde{\mathbf{B}}( Q)\rightarrow \Q(v^{\frac{1}{2}})$; see \eqref{eq:chi1}--\eqref{eq:chi3}. The formula of $\chi$ is extraordinary, which involves the polynomials $a_{2\lambda_f}(v),a_{\lambda_f}(v^2)$ to count the cardinality of automorphism groups of modules over $\F_q$ and $\F_{q^2}$, and also the square root $b_{\lambda_{c}}(v)$ of $a_{\lambda_c}(v)$.   
	
	In Theorem \ref{thm:coproduct}, an algebra homomorphism $\widetilde{\Delta}: \widetilde{\ch}(Q,\varrho)\rightarrow \widetilde{\ch}(Q,\varrho)\otimes\widetilde{\ch}(Q^{\dbl},\swa)$ is constructed such that the following diagram commutes:	\[ \xymatrix{ \tUi \ar[rr]^{\Delta\circ\imath }\ar[d]^{\widetilde{\psi}} && \tUi\otimes\tU \ar[d]^{\widetilde{\psi}\otimes\widetilde{\psi}^{\dbl}} 
		\\
		\widetilde{\ch}(Q,\varrho) \ar[rr]^{\widetilde{\Delta}\qquad\qquad} && \widetilde{\ch}(Q,\varrho)\otimes\widetilde{\ch}(Q^{\dbl},\swa)} \]
	This result interprets the coideal algebra structure of $\tUi$ in $\tU$. 
	
	It is striking that the maps $\widetilde{\Omega}$ and $\widetilde{\Delta}$ admit neat, compact closed-form expressions on the Hall bases. We expect that these maps can be formulated for arbitrary quivers. However, the $\varrho$-twisted compatible linear map $\chi$ in this general setting would be much more complicated than this one given in \eqref{eq:chi1}--\eqref{eq:chi3}. 
	
	We give an application for integral forms and dual canonical bases of $\tU$ and $\tUi$.
	Let $\cz=\Z[v^{\frac12},v^{-\frac12}]$. An integral form $\widetilde{\ch}(Q,\varrho)_\cz$ of $\widetilde{\ch}(Q,\varrho)$ is introduced in \cite{LP25}, and it has the dual canonical basis; similar for $\widetilde{\ch}(Q^{\dbl},\swa)$. Then it induces the integral form $\tUi_\cz$ of $\tUi$ (and also dual canonical basis) via the isomorphism $\widetilde{\psi}$. We prove that the homomorphisms $\widetilde{\Omega}: \widetilde{\ch}(Q,\varrho)\rightarrow \widetilde{\ch}(Q^{\dbl},\swa)$ and $\widetilde{\Delta}: \widetilde{\ch}(Q,\varrho)\rightarrow \widetilde{\ch}(Q,\varrho)\otimes\widetilde{\ch}(Q^{\dbl},\swa)$ preserve their integral forms in Proposition \ref{prop:interal-embedding-coproduct}. The key point is to prove that the $\varrho$-twisted compatible linear map $\chi:\widetilde{\mathbf{B}}( Q)\rightarrow \Q(v^{\frac{1}{2}})$ maps Hall basis elements to $\cz$; see Lemma \ref{lem:laurent-chi}. 
	
	Let $\mathbf{C}^\imath=\{\K_\alpha\diamond \mathfrak{L}_\lambda\mid \alpha\in\Z^\I,\lambda\in\mathfrak{P}\}$ (resp. $\mathbf{C}=\{\K_\alpha\diamond \K_{\beta^\diamond}\diamond \mathfrak{L}_{\lambda,\mu^\diamond}\mid \alpha,\beta\in\Z^\I,\lambda,\mu\in\fp\}$) be the dual canonical basis of $\tUi$ (resp. $\tU$). For
	\begin{align*}
		\imath(b)=\sum_{c\in\mathbf{C}}f_c(v)c,\qquad
		\widetilde{\Delta}(b)=  \sum_{b'\in\mathbf{C}^\imath,c\in\mathbf{C}}g_{b',c''}(v)b'\otimes c,
	\end{align*}
	we can see all the constants $f_c(v),g_{b',c''}(v)\in\cz$; see Proposition \ref{prop:integral-dCB}. 
	In fact, it is conjectured in \cite{LP25b} that these coefficients are in $\N[v^{\frac12},v^{-\frac12}]$. In order to solve this question, we shall develop the geometric realization of $\tU$ and $\tUi$ in \cite{Qin,LW21b} further to interpret the comultiplication and embedding $\imath:\tUi\to \tU$ by perverse sheaves; cf. \cite{Lus90,Lus93,Na01,FL21}.

	This paper is organized as follows. In Section \ref{sec:QG-iQG}, we review the basics of quantum groups and iquantum groups. In Section \ref{sec:iHall}, we review the construction of iquiver algebras and their iHall algebras, and the Hall algebra realization of quantum groups and iquantum groups. 
	
	In Section \ref{sec:main}, we give the Hopf algebra structure of Hall algebra $\widetilde{\ch}(Q,\varrho)$, realize the quantum symmetric pair $(\tU,\tUi)$ and interpret the coideal algebra structure of $\tUi$ in $\tU$ via Hall algebras. 
	In Section \ref{sec:dCB}, we consider integral forms and dual canonical bases as an application.

	\vspace{2mm}
	\noindent {\bf Acknowledgement.} ML deeply thanks Liangang Peng and Weiqiang Wang for their guidance and collaboration in related projects. ML is partially supported by the National Natural Science Foundation of China (No. 12171333).

	\section{Quantum groups and {\rm i}quantum groups}
	\label{sec:QG-iQG}
	
	\subsection{Quantum groups}
	\label{subsec:QG}

	Let $\I=\{1,\dots,n\}$ be the index set. 
	Let $C=(c_{ij})_{i,j \in \I}$ be the Cartan matrix of a simple Lie algebra $\fg$ of type ADE. Its Dynkin diagram is  denoted by $\Gamma$. 
	Let $\Delta=\{\alpha_i\mid i\in\I\}$ be the set of simple roots of $\fg$, and denote the root lattice by $\Z^{\I}:=\Z\alpha_1\oplus\cdots\oplus\Z\alpha_n$. We define a symmetric bilinear form on $\Z^\I$ by setting
	\begin{align} \label{BilForm}
		(\alpha_i,\alpha_j)=c_{ij},\quad \forall i,j\in\I.
	\end{align}

	Let $v$ be an indeterminate. 
	Denote, for $r\in\N,m \in \Z$,
	\[
	[r]:=[r]_{v}=\frac{v^r-v^{-r}}{v-v^{-1}},
	\;
	[r]!:=[r]_{v}^!=\prod_{i=1}^r [i]_{v}, \; \qbinom{m}{r}:=\qbinom{m}{r}_{v} =\frac{[m]_{v}[m-1]_{v}\ldots [m-r+1]_{v}}{[r]_{v}^!}.
	\]


	Following \cite{Dr87, BG17a}, the (Drinfeld double) quantum group $\tU := \tU_v(\fg)$ is defined to be the $\Q(v^{\frac{1}{2}})$-algebra generated by $E_i,F_i, \tK_i^{\pm1},\tK_i'^{\pm1}$, $i\in \I$, subject to the following relations:  for $i, j \in \I$,
	\begin{align}
		[E_i,F_j]= \delta_{ij}(v^{-1}-v) (\tK_i-\tK_i'),  &\qquad [\tK_i,\tK_j]=[\tK_i,\tK_j']  =[\tK_i',\tK_j']=0,
		\label{eq:KK}
		\\
		\tK_i E_j=v^{c_{ij}} E_j \tK_i, & \qquad \tK_i F_j=v^{-c_{ij}} F_j \tK_i,
		\label{eq:EK}
		\\
		\tK_i' E_j=v^{-c_{ij}} E_j \tK_i', & \qquad \tK_i' F_j=v^{c_{ij}} F_j \tK_i',
		\label{eq:K2}
	\end{align}
	and for $i\neq j \in \I$,
	\begin{align}
		& \sum_{r=0}^{1-c_{ij}} (-1)^r \left[ \begin{array}{c} 1-c_{ij} \\r \end{array} \right]_{v}  E_i^r E_j  E_i^{1-c_{ij}-r}=0,
		\label{eq:serre1} 
		\\ 
		&\sum_{r=0}^{1-c_{ij}} (-1)^r \left[ \begin{array}{c} 1-c_{ij} \\r \end{array} \right]_{v}  F_i^r F_j  F_i^{1-c_{ij}-r}=0.
		\label{eq:serre2}
	\end{align}
	
	The Drinfeld--Jimbo quantum group $\bU$ is defined to be the $\Q(v^{\frac12})$-algebra generated by $E_i,F_i, K_i, K_i^{-1}$, $i\in \I$, subject to the relations modified from \eqref{eq:KK}--\eqref{eq:serre2} with $\tK_i'$ replaced by $K_i^{-1}$. We can view $\bU$ as the quotient algebra of $\tU$ modulo the ideal $( \tK_i \tK_i'- 1 \mid  i\in \I )$; see \cite{Dr87}.

	Let $\widetilde{\bU}^+$ be the subalgebra of $\widetilde{\bU}$ generated by $E_i$ $(i\in \I)$, $\widetilde{\bU}^0$ be the subalgebra of $\widetilde{\bU}$ generated by $\tK_i, \tK_i'$ $(i\in \I)$, and $\widetilde{\bU}^-$ be the subalgebra of $\widetilde{\bU}$ generated by $F_i$ $(i\in \I)$, respectively.
	The subalgebras $\bU^+$, $\bU^0$ and $\bU^-$ of $\bU$ are defined similarly. Then both $\widetilde{\bU}$ and $\bU$ have triangular decompositions:
	$
	\widetilde{\bU} =\widetilde{\bU}^+\otimes \widetilde{\bU}^0\otimes\widetilde{\bU}^-,
	\, 
	\bU =\bU^+\otimes \bU^0\otimes\bU^-.
	$ 
	Clearly, ${\bU}^+\cong\widetilde{\bU}^+$, ${\bU}^-\cong \widetilde{\bU}^-$, and ${\bU}^0 \cong \widetilde{\bU}^0/(\tK_i \tK_i' -1 \mid   i\in \I)$.

	The algebras $\tU$ (and $\U$) are Hopf algebras, with the coproduct $\Delta$, the counit $\varepsilon$ and the antipode $S$ defined by
	\begin{align}
		\begin{split}
			\Delta(E_i)=E_i\otimes 1+K_i\otimes E_i,\quad &\Delta(F_i)=1\otimes F_i+F_i\otimes K_i',
			\\
			\Delta(K_i)=K_i\otimes K_i,\quad &\Delta(K_i')=K_i'\otimes K_i';
			\\
			\varepsilon(E_i)=0=\varepsilon(F_i),\qquad &\varepsilon(K_i)=1=\varepsilon(K_i');
			\\
			S(E_i)=-K_i^{-1}E_i,\quad S(F_i)=-F_iK_i'^{-1},\quad &S(K_i)=K_i^{-1},\quad S(K_i')=K_i'^{-1}.
		\end{split}
	\end{align}

	\begin{lemma}[\cite{BG17a}] \label{QG bar-involution def}
		
		There exists an anti-involution (called the bar-involution) $u\mapsto \ov{u}$ on $\tU$ (and $\U$) given by $\ov{v^{1/2}}=v^{-1/2}$, $\ov{E_i}=E_i$, $\ov{F_i}=F_i$, and $\ov{K_i}=K_i$, $\ov{K_i'}=K_i'$, for $i\in\I$.
	\end{lemma}

	\subsection{iQuantum groups}
	
	For a Cartan matrix $C=(c_{ij})_{i,j\in \I}$, let $\text{Inv}(C)$ be the group of permutations $\btau$ of the set $\I$ such that $c_{ij}=c_{\btau i,\btau j}$, for all $i,j$, and $\btau^2=\Id$. Then  $\btau \in \text{Inv}(C)$ can be viewed as an involution (which is allowed to be the identity) of the corresponding Dynkin diagram (which is identified with $\I$ by abuse of notation). We shall refer to the pair $(\I,\btau)$ as a (quasi-split) Satake diagram. 
	

	Associated with the Satake diagram $(\I,\varrho)$, following \cite{LW22a} we define 
	${\tU}^\imath$ to be the $\Q(v^{\frac12})$-subalgebra of $\tU$ generated by
	\begin{equation}
		\label{eq:Bi}
		B_i= F_i +  
		E_{\btau i} \tK_i',
		\qquad \tk_i = \tK_i \tK_{\btau i}', 
		\quad \forall i \in \I, 
	\end{equation}
	with $\tk_i$ invertible. Let $\bvs=(\vs_i)\in(\Q(v^{\frac12})^\times)^\I$ be such that $\vs_i=\vs_{\varrho i}$ for each $i\in\I$ which satisfies $c_{i,\varrho i}=0$. The iquantum groups \`a la Letzter-Kolb \cite{Let99,Ko14} $\Ui=\Ui_{\bvs}$ is the $\Q(v^{\frac12})$-subalgebra of $\U$ generated by
	$$B_i=F_i+\vs_i E_{\varrho i}K_i^{-1},\quad k_i=K_iK_{\btau i}^{-1},\quad \forall i\in\I.$$ 
	By \cite[Proposition 6.2]{LW22a}, the $\Q(v^{\frac12})$-algebra $\Ui$ is isomorphic to the quotient of $\tUi$ by the ideal generated by $\tk_i-\vs_i$ (for $i=\btau i$) and $\tk_i\tk_{\btau i}-\vs_i\vs_{\tau i}$ (for $i\neq \btau i$).
	
	Denote by $\imath:\tUi\to \tU$, $\imath:\Ui\to \U$ the natural embeddings.
	The algebras $\widetilde{\bU}^\imath$ and $\Ui$ are right coideal subalgebras of $\widetilde{\bU}$ and $\U$. respectively. 
	The pairs $(\widetilde{\bU}, \widetilde{\bU}^\imath)$ and $(\U,\Ui)$ are called quantum symmetric pairs, and $\tUi$ and $\Ui$ are called the universal {\em (quasi-split) iquantum groups}; they are {\em split} if $\btau =\Id$.

	For $i\in\I$, for any $\alpha=\sum_{i\in\I} a_i\alpha_i\in\Z^\I$, we set
	\begin{align}
		\label{eq:bbKi}
		\K_i=v^{\frac{1}{2}c_{i,\btau i}}\tk_i,\qquad \K_\alpha:=\prod_{i\in\I}\K_i^{a_i}.
	\end{align}

	\begin{lemma} [\text{cf. \cite[Lemma 3.6]{CLPRW}}]
		\label{lem:involution-iQG}
		There exists an anti-involution (called bar-involution) $:u\mapsto \ov{u}$ on $\tUi$ given by $\ov{v^{1/2}}=v^{-1/2}$, $\ov{B_i}=B_i$, and $\ov{\K_i}=\K_i$,  for $i\in\I$. In particular, $\ov{\tk_i}=v^{c_{i,\btau i}}\tk_{i}$. 
	\end{lemma}

	\begin{example} 
		Quasi-split iquantum groups of finite type are associated to the Satake diagrams of finite type, i.e., Dynkin graphs with involution. The split Satake diagrams formally look the same as Dynkin diagrams. A complete list of such connected quasi-split Satake diagrams of finite type with $\btau \neq \Id$ is given in Table \ref{tab:Satakediag} below.  
		\begin{table}[h]  
			\begin{center}
				\centering
				\begin{tabular}{|m{5cm}<{\centering}|m{8cm}<{\centering}|}
					\hline
					Types &  Satake diagrams  \\
					\hline

					${\rm AIII}_{2n-1}$ $(n\geq 1)$ & \setlength{\unitlength}{0.6mm}			\begin{picture}(70,35)(0,5)
						
						\put(0,10){$\circ$}
						\put(0,30){$\circ$}
						\put(50,10){$\circ$}
						\put(50,30){$\circ$}
						\put(72,10){$\circ$}
						\put(72,30){$\circ$}
						\put(92,20){$\circ$}
						\put(-5,5.5){$2n-1$}
						\put(-2,34){${1}$}
						\put(47,6){\small $n+2$}
						\put(47,34){\small $n-2$}
						\put(69,6){\small $n+1$}
						\put(69,34){\small $n-1$}
						\put(92,16){\small $n$}
						
						\put(3,11.5){\line(1,0){16}}
						\put(3,31.5){\line(1,0){16}}
						\put(23,10){$\cdots$}
						\put(23,30){$\cdots$}
						\put(33.5,11.5){\line(1,0){16}}
						\put(33.5,31.5){\line(1,0){16}}
						\put(53,11.5){\line(1,0){18.5}}
						\put(53,31.5){\line(1,0){18.5}}
						
						\put(75,12){\line(2,1){17}}
						\put(75,31){\line(2,-1){17}}

						\color{red}
						\put(-7,20){$\btau$}
						\qbezier(0,13.5)(-4,21.5)(0,29.5)
						\put(-0.25,14){\vector(1,-2){0.5}}
						\put(-0.25,29){\vector(1,2){0.5}}
						
						\qbezier(50,13.5)(46,21.5)(50,29.5)
						\put(49.75,14){\vector(1,-2){0.5}}
						\put(49.75,29){\vector(1,2){0.5}}

						\qbezier(72,13.5)(68,21.5)(72,29.5)
						\put(71.75,14){\vector(1,-2){0.5}}
						\put(71.75,29){\vector(1,2){0.5}}
					\end{picture}\\
					\hline
					${\rm AIII}_{2n}$ $(n\geq 1)$ & \setlength{\unitlength}{0.6mm}			\begin{picture}(70,35)(0,5)
						
						\put(0,10){$\circ$}
						\put(0,30){$\circ$}
						\put(50,10){$\circ$}
						\put(50,30){$\circ$}
						\put(72,10){$\circ$}
						\put(72,30){$\circ$}
						\put(-4,5){$2n$}
						\put(-2,34){${1}$}
						\put(47,6){\small $n+2$}
						\put(47,34){\small $n-1$}
						\put(69,6){\small $n+1$}
						\put(69,34){\small $n$}
						
						\put(3,11.5){\line(1,0){16}}
						\put(3,31.5){\line(1,0){16}}
						\put(23,10){$\cdots$}
						\put(23,30){$\cdots$}
						\put(33.5,11.5){\line(1,0){16}}
						\put(33.5,31.5){\line(1,0){16}}
						\put(53,11.5){\line(1,0){18.5}}
						\put(53,31.5){\line(1,0){18.5}}
						
						\put(73.5,13.6){\line(0,1){16}}
						
						\color{red}
						\put(-7,20){$\btau$}
						\qbezier(0,13.5)(-4,21.5)(0,29.5)
						\put(-0.25,14){\vector(1,-2){0.5}}
						\put(-0.25,29){\vector(1,2){0.5}}
						
						\qbezier(50,13.5)(46,21.5)(50,29.5)
						\put(49.75,14){\vector(1,-2){0.5}}
						\put(49.75,29){\vector(1,2){0.5}}

						\qbezier(76,13.5)(80,21.5)(76,29.5)
						\put(75.85,14){\vector(-1,-2){0.5}}
						\put(75.85,29){\vector(-1,2){0.5}}
					\end{picture}\\
					\hline
					${\rm DI}_n$ $(n\geq4)$ & \setlength{\unitlength}{0.6mm}	\begin{picture}(40,35)(20,-15)
						\put(0,-1){$\circ$}
						\put(0,-6){\small$1$}
						
						\put(3,0){\line(1,0){16.5}}
						\put(20,-1){$\circ$}
						\put(20,-6){\small$2$}
						\put(64,-1){$\circ$}
						\put(56,-6){\small$n-2$}
						\put(84,-10){$\circ$}
						\put(80,-13){\small${n-1}$}
						\put(84,9.5){$\circ$}
						\put(84,13.5){\small${n}$}

						\put(38,0){\line(-1,0){15.5}}
						\put(64,0){\line(-1,0){15}}
						
						\put(40,-1){$\cdots$}
						
						\put(83.5,9.5){\line(-2,-1){16.5}}
						\put(83.5,-8.5){\line(-2,1){16.5}}
						
						\put(12,-20.5){\begin{picture}(100,40)	\color{red}
								\put(79,20){$\btau$}
								\qbezier(75,13.5)(79,21.5)(75,29.5)
								\put(75.25,14){\vector(-1,-2){0.5}}
								\put(75.25,29){\vector(-1,2){0.5}}
						\end{picture}}
					\end{picture}
					\\
					\hline
					${\rm EII}_6$  &	\setlength{\unitlength}{0.6mm}		\begin{picture}(70,35)(20,7)
						\put(97,36){\small${6}$}
						\put(75,36){\small${5}$}
						\put(97,6.5){\small${1}$}
						\put(75,6.5){\small${2}$}
						\put(33,16){\small $4$}	
						\put(55,16){\small $3$}	\put(10,35){\rotatebox[origin=c]{180}{\begin{picture}(100,10)
									
									\put(10,10){$\circ$}
									
									\put(32,10){$\circ$}
									
									\put(10,30){$\circ$}
									
									\put(32,30){$\circ$}

									
									\put(31.5,11){\line(-1,0){19}}
									\put(31.5,31){\line(-1,0){19}}
									
									\put(52,22){\line(-2,1){17.5}}
									\put(52,20){\line(-2,-1){17.5}}
									
									\put(54.7,21.2){\line(1,0){19}}

									\put(52,20){$\circ$}
									
									\put(74,20){$\circ$}
							\end{picture}}

							\put(-37,-32.5){\begin{picture}(100,40)	\color{red}
									\qbezier(5,13.5)(9,21.5)(5,29.5)
									\put(5.25,14){\vector(-1,-2){0.5}}
									\put(5.25,29){\vector(-1,2){0.5}}
							\end{picture}}
							\put(-15,-32.5){\begin{picture}(100,40)	\color{red}
									\put(9,20){$\btau$}
									\qbezier(5,13.5)(9,21.5)(5,29.5)
									\put(5.25,14){\vector(-1,-2){0.5}}
									\put(5.25,29){\vector(-1,2){0.5}}
							\end{picture}}
						}
						
					\end{picture}
					\\\hline
				\end{tabular} 
			\end{center}
			\vspace{0.5cm}
			\caption{Quasi-split Satake diagrams of finite type with $\btau \neq {\rm Id}$}
			\label{tab:Satakediag}
		\end{table}
	\end{example}
	
	\begin{example}
		\label{ex:QGvsiQG}
		{\rm (Quantum groups as iquantum groups of diagonal type)} 
		Consider the $\Q(v^{\frac12})$-subalgebra $\tUUi$ of $\tUU$
		generated by
		\[
		\ck_i:=\tK_{i} \tK_{i^{\diamond}}', \quad
		\ck_i':=\tK_{i^{\diamond}} \tK_{i}',  \quad
		\cb_{i}:= F_{i}+ E_{i^{\diamond}} \tK_{i}', \quad
		\cb_{i^{\diamond}}:=F_{i^{\diamond}}+ E_{i} \tK_{i^{\diamond}}',
		\qquad \forall i\in \I.
		\]
		Here we drop the tensor product notation and use instead $i^\diamond$ to index the generators of the second copy of $\tU$ in $\tUU$. There exists a $\Q(v^{\frac12})$-algebra isomorphism $\widetilde{\phi}: \tU \rightarrow \tUUi$ such that
		\[
		\widetilde{\phi}(E_i)= \cb_{i},\quad \widetilde{\phi}(F_i)= \cb_{i^{\diamond}}, \quad \widetilde{\phi}(\tK_i)= \ck_i', \quad \widetilde{\phi}(\tK_i')= \ck_i, \qquad \forall  i\in \I.
		\]
		In this case, the Satake diagram is $(\Gamma\sqcup\Gamma^\diamond,\swa)$, where $\Gamma^\diamond$ is a copy of the Dynkin diagram $\Gamma$ of $\tU$. 
	\end{example}

	\section{{\rm i}Quiver algebras and {\rm i}Hall algebras}
	\label{sec:iHall}
	
	In this section, we recall the construction of iquiver algebras and the realization of iquantum groups via iHall algebras formulated in \cite{LW22a}. Throughout this paper, let $\bfk=\mathbb{F}_q$ be a finite field of $q$ elements.

	\subsection{The iquivers and doubles}
	\label{subsec:i-quiver}
	
	Let $Q=(Q_0,Q_1)$ be a Dynkin quiver. 
	An {\em involution} of $Q$ is defined to be an automorphism $\btau$ of the quiver $Q$ such that $\btau^2=\Id$ (we allow the {\em trivial} involution $\Id:Q\rightarrow Q$). An involution $\btau$ of $Q$ induces an involution of the path algebra $\bfk Q$, again denoted by $\btau$.
	A quiver together with an involution $\btau$, $(Q, \btau)$, will be called an {\em iquiver}. 
	
	Let $R_2$ denote the radical square zero of the path algebra of $\xymatrix{1 \ar@<0.5ex>[r]^{\varepsilon} & 1' \ar@<0.5ex>[l]^{\varepsilon'}}$, 
	i.e., $\varepsilon' \varepsilon =0 =\varepsilon\varepsilon '$. Define a $\bfk$-algebra
	\begin{equation}
		\label{eq:La}
		\Lambda=\bfk Q\otimes_\bfk R_2.
	\end{equation}
	
	Associated to the quiver $Q$, the {\em double framed quiver} $Q^\sharp$
	is the quiver such that
	\begin{itemize}
		\item[(i)] the vertex set of $Q^{\sharp}$ consists of 2 copies of the vertex set $Q_0$, $\{i,i'\mid i\in Q_0\}$;
		\item[(ii)] the arrow set of $Q^{\sharp}$ is
		\[
		\{\alpha: i\rightarrow j,\alpha': i'\rightarrow j'\mid(\alpha:i\rightarrow j)\in Q_1\}\cup\{ \varepsilon_i: i\rightarrow i' ,\varepsilon'_i: i'\rightarrow i\mid i\in Q_0 \}.
		\]
	\end{itemize}
	Let $I^{\sharp}$ be the admissible ideal of $\bfk Q^{\sharp}$ generated by
	\begin{itemize}
		\item[(1)]
		(Nilpotent relations) $\varepsilon_i \varepsilon_i'$, $\varepsilon_i'\varepsilon_i$ for any $i\in Q_0$;
		\item[(2)]
		(Commutative relations) $\varepsilon_j' \alpha' -\alpha\varepsilon_i'$, $\varepsilon_j \alpha -\alpha'\varepsilon_i$ for any $(\alpha:i\rightarrow j)\in Q_1$.
	\end{itemize}
	Then the algebra $\La$ can be realized as $\Lambda\cong \bfk Q^{\sharp} \big/ I^{\sharp}$.
	
	Note $Q^\sharp$ admits a natural involution, $\swa$.
	The involution $\btau$ of a quiver $Q$ induces an involution ${\btau}^{\sharp}$ of $Q^{\sharp}$ which is the composition of $\swa$ and $\btau$ (on the two copies of subquivers $Q$ and $Q'$ of $Q^\sharp$).
	By \cite[Definition 2.5]{LW22a}, the {\rm iquiver algebra} of $(Q, \btau)$ is the fixed point subalgebra of $\Lambda$ under ${\btau}^{\sharp}$,
	\begin{equation}
		\label{eq:iLa}
		\iLa
		= \{x\in \Lambda\mid {\btau}^{\sharp}(x) =x\}.
	\end{equation}
	The algebra $\iLa$ can be described in terms of a certain quiver $\ov{Q}$ and its admissible ideal $\ov{I}$ so that $\iLa \cong \bfk \ov{Q} / \ov{I}$; see \cite[Proposition 2.6]{LW22a}.
	We recall $\ov{Q}$ and $\ov{I}$ as follows:
	\begin{itemize}
		\item[(i)] $\ov{Q}$ is constructed from $Q$ by adding a loop $\varepsilon_i$ at the vertex $i\in Q_0$ if $\btau i=i$, and adding an arrow $\varepsilon_i: i\rightarrow \btau i$ for each $i\in Q_0$ if $\btau i\neq i$;
		\item[(ii)] $\ov{I}$ is generated by
		\begin{itemize}
			\item[(1)] (Nilpotent relations) $\varepsilon_{i}\varepsilon_{\btau i}$ for any $i\in\I$;
			\item[(2)] (Commutative relations) $\varepsilon_i\alpha-\btau(\alpha)\varepsilon_j$ for any arrow $\alpha:j\rightarrow i$ in $Q_1$.
		\end{itemize}
	\end{itemize}

	The following quivers are examples of the quivers $\ov{Q}$ used to describe the iquiver algebras $\Lambda^\imath$
	associated to non-split iquivers of type ADE; cf. \cite{LW22a}.

	\begin{center}\setlength{\unitlength}{0.7mm}
		\vspace{-1.5cm}
		\begin{equation}
			\label{diag: A}
			\begin{picture}(100,40)(0,20)
				\put(0,10){$\circ$}
				\put(0,30){$\circ$}
				\put(50,10){$\circ$}
				\put(50,30){$\circ$}
				\put(72,10){$\circ$}
				\put(72,30){$\circ$}
				\put(92,20){$\circ$}
				\put(-8,6){$2n+1$}
				\put(-0,34){${1}$}
				\put(45,6){\small $n+3$}
				\put(45,34){\small ${n-1}$}
				\put(67,6){\small $n+2$}
				\put(72,34){\small ${n}$}
				\put(92,16){\small $n+1$}
				
				\put(3,11.5){\vector(1,0){16}}
				\put(3,31.5){\vector(1,0){16}}
				\put(23,10){$\cdots$}
				\put(23,30){$\cdots$}
				\put(33.5,11.5){\vector(1,0){16}}
				\put(33.5,31.5){\vector(1,0){16}}
				\put(53,11.5){\vector(1,0){18.5}}
				\put(53,31.5){\vector(1,0){18.5}}
				\put(75,12){\vector(2,1){17}}
				\put(75,31){\vector(2,-1){17}}
				\color{purple}
				\put(0,29.5){\vector(0,-1){17}}
				\put(2,13){\vector(0,1){17}}
				\put(50,29.5){\vector(0,-1){17}}
				\put(52,13){\vector(0,1){17}}
				\put(72,29.5){\vector(0,-1){17}}
				\put(74,13){\vector(0,1){17}}
				
				\put(-5,20){$\varepsilon_{1}$}
				\put(3,20){$\varepsilon_{2n+1}$}
				\put(39,20){\small $\varepsilon_{n-1}$}
				\put(53,20){\small $\varepsilon_{n+3}$}
				\put(66,20){\small $\varepsilon_{n}$}
				\put(75,20){\small $\varepsilon_{n+2}$}
				\put(92,30){\small $\varepsilon_{n+1}$}
				
				\qbezier(93,23)(90.5,25)(92,27)
				\qbezier(92,27)(94,30)(97,27)
				\qbezier(97,27)(98,24)(95.5,22.6)
				\put(95.6,23){\vector(-1,-1){0.3}}
			\end{picture}
		\end{equation}
		\vspace{-0.6cm}
	\end{center}

	\begin{center}\setlength{\unitlength}{0.8mm}
		\begin{equation}
			\label{diag: D}
			\begin{picture}(100,25)(-5,0)
				\put(0,-1){$\circ$}
				\put(0,-5){\small$1$}
				\put(20,-1){$\circ$}
				\put(20,-5){\small$2$}
				\put(64,-1){$\circ$}
				\put(84,-10){$\circ$}
				\put(80,-13){\small${n-1}$}
				\put(84,9.5){$\circ$}
				\put(84,12.5){\small${n}$}

				\put(19.5,0){\vector(-1,0){16.8}}
				\put(38,0){\vector(-1,0){15.5}}
				\put(64,0){\vector(-1,0){15}}
				
				\put(40,-1){$\cdots$}
				\put(83.5,9.5){\vector(-2,-1){16}}
				\put(83.5,-8.5){\vector(-2,1){16}}
				\color{purple}
				\put(86,-7){\vector(0,1){16.5}}
				\put(84,9){\vector(0,-1){16.5}}
				
				\qbezier(63,1)(60.5,3)(62,5.5)
				\qbezier(62,5.5)(64.5,9)(67.5,5.5)
				\qbezier(67.5,5.5)(68.5,3)(66.4,1)
				\put(66.5,1.4){\vector(-1,-1){0.3}}
				\qbezier(-1,1)(-3,3)(-2,5.5)
				\qbezier(-2,5.5)(1,9)(4,5.5)
				\qbezier(4,5.5)(5,3)(3,1)
				\put(3.1,1.4){\vector(-1,-1){0.3}}
				\qbezier(19,1)(17,3)(18,5.5)
				\qbezier(18,5.5)(21,9)(24,5.5)
				\qbezier(24,5.5)(25,3)(23,1)
				\put(23.1,1.4){\vector(-1,-1){0.3}}
				
				\put(-1,9.5){$\varepsilon_1$}
				\put(19,9.5){$\varepsilon_2$}
				\put(59,9.5){$\varepsilon_{n-2}$}
				\put(79,-1){$\varepsilon_{n}$}
				\put(87,-1){$\varepsilon_{n-1}$}
			\end{picture}
		\end{equation}
		\vspace{.8cm}
	\end{center}

	\begin{center}\setlength{\unitlength}{0.8mm}
		\vspace{-3cm}
		\begin{equation}
			\label{diag: E}
			\begin{picture}(100,40)(0,20)
				\put(10,6){\small${6}$}
				\put(10,10){$\circ$}
				\put(32,6){\small${5}$}
				\put(32,10){$\circ$}
				
				\put(10,30){$\circ$}
				\put(10,33){\small${1}$}
				\put(32,30){$\circ$}
				\put(32,33){\small${2}$}
				
				\put(31.5,11){\vector(-1,0){19}}
				\put(31.5,31){\vector(-1,0){19}}
				
				\put(52,22){\vector(-2,1){17.5}}
				\put(52,20){\vector(-2,-1){17.5}}
				
				\put(54.7,21.2){\vector(1,0){19}}
				
				\put(52,20){$\circ$}
				\put(52,16.5){\small$3$}
				\put(74,20){$\circ$}
				\put(74,16.5){\small$4$}
				\color{purple}
				\put(10,12.5){\vector(0,1){17}}
				\put(12,29.5){\vector(0,-1){17}}
				\put(32,12.5){\vector(0,1){17}}
				\put(34,29.5){\vector(0,-1){17}}
				
				\qbezier(52,22.5)(50,24)(51,26.5)
				\qbezier(51,26.5)(53,29)(56,26.5)
				\qbezier(56,26.5)(57.5,24)(55,22)
				\put(55.1,22.4){\vector(-1,-1){0.3}}
				\qbezier(74,22.5)(72,24)(73,26.5)
				\qbezier(73,26.5)(75,29)(78,26.5)
				\qbezier(78,26.5)(79,24)(77,22)
				\put(77.1,22.4){\vector(-1,-1){0.3}}
				
				\put(35,20){$\varepsilon_2$}
				\put(27,20){$\varepsilon_5$}
				\put(13,20){$\varepsilon_1$}
				\put(5,20){$\varepsilon_6$}
				\put(52,30){$\varepsilon_3$}
				\put(73,30){$\varepsilon_4$}
			\end{picture}
		\end{equation}
		\vspace{1cm}
	\end{center}

	
	\begin{example}[iQuivers of diagonal type]
		\label{ex:diagquiver}
		Let $Q$ be an arbitrary quiver, and $Q^{\dbl} =Q\sqcup  Q^{\diamond}$,  where $Q^{\diamond}$ is an identical copy of $Q$  with a vertex set $\{i^{\diamond} \mid i\in Q_0\}$ and an arrow set $\{ \alpha^{\diamond} \mid \alpha \in Q_1\}$. We let $\rm{swap}$ be the involution of $Q^{\rm dbl}$ uniquely determined by $\swa(i)=i^\diamond$ for any $i\in Q_0$. 
		Then $(Q^{\rm dbl},\mathrm{swap})$ is an iquiver with $\Lambda$  as its iquiver algebra; see \cite[Example 2.10]{LW22a}. 
	\end{example}

	By \cite[Corollary 2.12]{LW22a}, $\bfk Q$ is naturally a subalgebra and also a quotient algebra of $\Lambda^\imath$.
	Viewing $\bfk Q$ as a subalgebra of $\Lambda^{\imath}$, we have a restriction functor
	\[
	\res: \mod(\Lambda^{\imath})\longrightarrow\mod(\bfk Q).
	\]
	Viewing $\bfk Q$ as a quotient algebra of $\Lambda^{\imath}$, we obtain a pullback functor
	\begin{equation}\label{eqn:rigt adjoint}
		\iota:\mod(\bfk Q)\longrightarrow\mod(\Lambda^{\imath}).
	\end{equation}
	Hence a simple module $S_i (i\in Q_0)$ of $\bfk Q$ is naturally a simple $\iLa$-module.

	For each $i\in Q_0$, define a $\bfk$-algebra (which can be viewed as a subalgebra of $\iLa$)
	\begin{align}\label{dfn:Hi}
		\BH _i:=\left\{ \begin{array}{cc}  \bfk[\varepsilon_i]/(\varepsilon_i^2) & \text{ if }\btau i=i,
			\\
			\bfk(\xymatrix{i \ar@<0.5ex>[r]^{\varepsilon_i} & \btau i \ar@<0.5ex>[l]^{\varepsilon_{\btau i}}})/( \varepsilon_i\varepsilon_{\btau i},\varepsilon_{\btau i}\varepsilon_i)  &\text{ if } \btau i \neq i .\end{array}\right.
	\end{align}
	For $i\in \I$, define the indecomposable module over $\BH _i$
	\begin{align}
		\label{eq:E}
		\E_i =\begin{cases}
			\bfk[\varepsilon_i]/(\varepsilon_i^2), & \text{ if }\btau i=i;
			\\
			\xymatrix{\bfk\ar@<0.5ex>[r]^1 & \bfk\ar@<0.5ex>[l]^0} \text{ on the quiver } \xymatrix{i\ar@<0.5ex>[r]^{\varepsilon_i} & \btau i\ar@<0.5ex>[l]^{\varepsilon_{\btau i}} }, & \text{ if } \btau i\neq i.
		\end{cases}
	\end{align}
	Then $\E_i$, for $i\in Q_0$, can be viewed as a $\iLa$-module and will be called a {\em generalized simple} $\iLa$-module.
	
	\begin{lemma}[\text{\cite[Proposition 3.9]{LW22a}}]
		\label{cor: res proj}
		For any $M\in\mod(\Lambda^{\imath})$ the following are equivalent: (i) $\pd M<\infty$; 
		(ii) $\ind M<\infty$;
		(iii) $\pd M\leq1$;
		(iv) $\ind M\leq1$;
		(v) $\res_\BH (M)$ is projective as an $\BH $-module.
	\end{lemma}
	
	
	The algebra $\Lambda^\imath$ is $1$-Gorenstein; see \cite[Proposition 3.9]{LW22a}. 
	We denote by $\cp^{\leq 1}(\Lambda^\imath)$ the subcategory of $\Lambda^\imath$-modules with finite projective dimension (equivalently, with projective dimension at most $1$). We see $\K_i\in\cp^{\leq1}(\Lambda^\imath)$.

	\subsection{Hall algebras}
	
	Let $\ca$ be an essentially small exact category, linear over the finite field $\bfk=\F_q$.
	Assume that $\ca$ has finite morphism and extension spaces:
	$$|\Hom_\ca(A,B)|<\infty,\quad |\Ext^1_\ca(A,B)|<\infty,\,\,\forall A,B\in\ca.$$
	

	We denote by $\Iso(\ca)$ the set of isoclasses of objects of $\ca$.	Given objects $X,Y,Z\in\ca$, define $\Ext^1_\ca(X,Z)_Y\subseteq \Ext^1_\ca(X,Z)$ to be the subset parameterising extensions with the middle term isomorphic to $Y$. We define the Ringel-Hall algebra (also called Hall algebra) $\ch(\ca)$ to be the $\Q$-vector space whose basis is formed by the isomorphism classes $[X]$ of objects $X$ of $\ca$, with the multiplication
	defined by
	\begin{align}
		\label{eq:mult}
		[X]\bullet [Z]=\sum_{[Y]\in \Iso(\ca)}G_{XZ}^Y[Y],\text{ where }G_{XZ}^Y:=\frac{|\Ext_\ca^1(X,Z)_Y|}{|\Hom_\ca(X,Z)|}.
	\end{align}
	It is well known that
	the algebra $\ch(\ca)$ is associative and unital. The unit is given by $[0]$, where $0$ is the zero object of $\ca$; see \cite{Rin90,Br13}. 

	For any three objects $X,Y,Z$, let
	\begin{align}
		\label{eq:Fxyz}
		F_{XZ}^Y&:= \big |\{L\subseteq Y \mid L \cong Z,  Y/L\cong X\} \big |.
	\end{align}
	Let $\Aut(X)$ be the automorphism group of $X\in\ca$, and set $a_X=|\Aut(X)|$.  The Riedtmann-Peng formula states that
	\[
	F_{XZ}^Y= G_{XZ}^Y \cdot \frac{a_Y}{a_X a_Z}.
	\]
	
	For any $M_1,\dots,M_n$ and $M\in\ca$, we define $G_{M_1,\dots,M_n}^M$ and $F_{M_1,\dots,M_n}^M$ recursively by
	\begin{align}
		\label{eq:iterated-Hall}
		G_{M_1,\dots,M_n}^M:=\sum_{[L]} G_{M_1,L}^MG_{M_2,\dots,M_n}^L,\qquad F_{M_1,\dots,M_n}^M:=\sum_{[L]} F_{M_1,L}^MF_{M_2,\dots,M_n}^L.
	\end{align}
	
	\subsection{Hall algebras and iHall algebras}
	\label{subsec:iHall}
	
	For a (Dynkin) quiver $Q$, we denote by $\langle\cdot ,\cdot\rangle_Q$ the Euler form of $\mod(\bfk Q)$. Define
	\begin{align*}
		(x,y)_Q&= \langle x,y\rangle_Q+\langle y,x\rangle_Q.
	\end{align*}
	In particular, $(S_i,S_j)_Q=c_{ij}$ for any $i,j\in \I$, the entries of the Cartan matrix $C$. We identify the Grothendieck group $K_0(\mod(\bfk Q))$ with the root lattice $\Z^\I$, and denote by $d^M=(d^M_i)_{i \in \I }\in\Z^\I$ the corresponding class in $K_0(\mod(\bfk Q))$ for any $\bfk Q$-module $M$.
	
	Let $\ch(\Lambda^\imath)$ be the  Hall algebra of $\mod(\Lambda^\imath)$. Let $\sqq$ be a fixed square root of $q$. We define the twisted Hall algebra $\widetilde{\ch}(\Lambda^\imath)$ to be the $\Q(\sqq^{1/2})$-algebra on the same vector space as $\ch(\Lambda^\imath)$ with twisted multiplication given by
	$$[M]*[N]=\sqq^{\langle \res M,\res N\rangle_{Q}}[M]\bullet [N].$$
	
	Let $\mathcal{I}$ be the ideal of $\widetilde{\ch}(\Lambda^\imath)$ generated by all differences 
	$[L]-[K\oplus M]$ if there is a short exact sequence
	\begin{equation}
		\label{eq:ideal}
		0 \longrightarrow K \longrightarrow L \longrightarrow M \longrightarrow 0
	\end{equation}
	with $K\in \cp^{\leq 1}(\Lambda^\imath)$.
	Let	\begin{equation}
		\label{eq:Sca}
		\cs := \{ a[K] \in \widetilde{\ch}(\Lambda^\imath)/\mathcal{I} \mid a\in \Q(\sqq^{1/2})^\times, \pd K\leq1\}.
	\end{equation}
	By \cite{LP21,LW22a}, the right localization of
	$\ch(\Lambda^\imath)/\mathcal{I}$ with respect to $\cs$ exists, and will be denoted by $\widetilde{\ch}(\bfk Q,\btau)$ or $\cs\cd\widetilde{\ch}(\Lambda^\imath)$, called the iHall algebra (also called the twisted semi-derived Ringel-Hall algebra) of $\mod(\Lambda^\imath)$.

	
	We have $[\K_i]*[\K_j]=[\K_j]*[\K_i]=[\K_i\oplus \K_j]$ in $\widetilde{\ch}(\bfk Q,\varrho)$ for any $i,j\in\I$; cf. \cite[Lemma 4.7]{LW22a}. 
	For any $\alpha=\sum_{i\in\I}a_i\alpha_i\in\N^{\I}$, we define in $\widetilde{\ch}(\bfk Q,\btau)$: 
	$$[\K_\alpha]=[\oplus_{i\in\I}\K_i^{\oplus a_i}]=\prod_{i\in\I}[\K_i]^{a_i}.$$
	Similarly, one can define $[\K_\alpha]$ in $\widetilde{\ch}(\bfk Q,\btau)$ for $\alpha\in\Z^\I$.

	Moreover, we have
	\begin{align}
		\label{eq:KX=XK}
		[\K_\alpha]*[X]=\sqq^{( \btau\alpha, \widehat{X}  )_Q-( \alpha, \widehat{X}  )_Q}[X]*[\K_\alpha],\quad \forall X\in\mod(\bfk Q), \alpha\in\N^{\I}. 
	\end{align}

	\begin{lemma}[\text{cf. \cite[Proposition 4.9]{LW22a}}]
		\label{basis-iHall}
		The algebra $\widetilde{\ch}(\bfk Q,\btau)$ has a (Hall) basis given by
		\begin{align}
			\label{eq:Hallbasis-hat}
			\{[X]*[\K_\alpha]\mid X\in\mod(\bfk Q)\subseteq \mod(\Lambda^\imath), \alpha\in\Z^{\I}\}.
		\end{align}
		
	\end{lemma}

	Let $\widetilde{\ct}(\bfk Q,\btau)$ be the subalgebra of $\widetilde{\ch}(\bfk Q,\btau)$ generated by $[\K_\alpha]$, $\alpha\in K_0(\mod(\bfk Q))$, which is a Laurent polynomial algebra in $[\E_i]$, for $i\in \I$.

	\begin{proposition}[\text{\cite[Proposition 3.3]{LR24}}]
		\label{prop:iHallmult}
		For any $A,B\in\mod(\bfk Q)\subseteq \mod(\Lambda^\imath)$, we have in $\widetilde{\ch}(\bfk Q,\btau)$ 
		\begin{align*}
			[A]*[B]=&
			\sum_{[L],[M],[N],[X]} \sqq^{\langle X,M\rangle-\langle \varrho X,M\rangle-\langle A,B\rangle} q^{\langle N,L\rangle} F_{N,L}^M F_{X,N}^AF_{ L,\varrho X}^{B} 
			\frac{a_L a_N a_X}{a_M} \cdot[K_{\widehat{X}}]*[M]
		\end{align*}
		where the sum is over $ [L],[M],[N],[X]\in\Iso(\mod(\bfk Q))$.
	\end{proposition}

	
	\begin{lemma}[cf. \text{\cite[Theorem 7.7]{LW22a}}]
		\label{lem:Hall-iQG}
		Let $(Q, \btau)$ be a Dynkin iquiver. Then we have the following isomorphism $\widetilde{\psi}:\tUi|_{v=\sqq}\stackrel{\simeq}{\rightarrow} \widetilde{\ch}(\bfk Q,\btau)$ of $\Q({\sqq^{1/2}})$-algebras, which sends
		\begin{align}
			\label{eq:psi}
			B_i \mapsto \sqq^{-\frac{1}{2}}[S_i], \qquad 
			\tilde{k}_i \mapsto
			\begin{cases}
				[\K_i]&\text{if $\varrho i\neq i$},\\
				\sqq^{-1}[\K_{i}]&\text{if $\varrho i=i$}.
			\end{cases}
		\end{align}
	\end{lemma}
	In this way, we have $\widetilde{\psi}(\K_i)=[\K_i]$, for $i\in\I$, which convinces the notation.

	As Example \ref{ex:QGvsiQG} shows, we can view quantum group $\tU$ to be an iquantum group; correspondingly, 
	as Example \ref{ex:diagquiver} shows, we can view the algebra $\Lambda$ in \eqref{eq:La} to be the iquiver algebra of $(Q^{\rm dbl},\swa)$.
	From Lemma \ref{lem:Hall-iQG}, we have the following result.
	\begin{lemma}[Bridgeland's Theorem reformulated] 
		\label{lem:Bridgeland}
		Let $Q$ be a Dynkin quiver. Then we have the following isomorphism of $\Q(\sqq^{1/2})$-algebras
		\begin{align*}
			&\widetilde{\psi}^{\dbl}:\tU|_{v=\sqq}\stackrel{\simeq}{\longrightarrow} \widetilde{\ch}(\bfk Q^{\rm dbl},\swa),
			\\
			E_i \mapsto \sqq^{-\frac{1}{2}}[S_i],&\qquad F_i\mapsto \sqq^{-\frac{1}{2}}[S_{i^\diamond}],
			\qquad
			K_i\mapsto [\K_{i^\diamond}],\qquad K_i'\mapsto [\K_i].
		\end{align*}
	\end{lemma}
	

	
	
	\subsection{Generic iHall algebras}
	\label{sub:generic}
	
	For a Dynkin iquiver $(Q,\btau)$,
	we recall the generic iHall algebras defined in \cite[\S9.3]{LW22a}.
	Recall that $\Phi^+$ is the set of positive roots.
	For any $\beta\in\Phi^+$, denote by $M_q(\beta)$ its corresponding indecomposable $\bfk Q$-module, i.e., $\dimv M_q(\beta)=\beta$.
	Let $\mathfrak{P}:=\mathfrak{P}(Q)$ be the set of functions $\lambda: \Phi^+\rightarrow \N$.
	Then the modules
	\begin{align}
		\label{def:Mlambda}
		M_q(\lambda):= \bigoplus_{\beta\in\Phi^+} M_q(\beta)^{\oplus \lambda(\beta)},\quad \text{ for } \lambda\in\mathfrak{P},
	\end{align}
	provide a complete set of isoclasses of $\bfk Q$-modules. We denote $d^\lambda=(d^\lambda_i)_{i \in \I} \in \N^\I$ the dimension vector of $M_q(\lambda)$. Similarly, for $\beta \in \Z^\I$ we also denote $\beta=(d^{\beta}_i)_{i \in \I}$.
	
	For any $\lambda,\mu,\nu\in\fp$, by \cite{Rin90}, there exist polynomials $f_{\lambda,\mu}^\nu(v),g_{\lambda,\mu}^\nu(v),a_\lambda(v)\in\Z[v,v^{-1}]$ such that 
	\begin{align*}
		F_{M_q(\lambda),M_q(\mu)}^{M_q(\nu)}=f_{\lambda,\mu}^\nu(\sqq),\quad G_{M_q(\lambda),M_q(\mu)}^{M_q(\nu)}=g_{\lambda,\mu}^\nu(\sqq),\quad a_{M_q(\lambda)}=a_\lambda(\sqq).
	\end{align*}
	In particular,
	we have
	\begin{align*}
		f_{\lambda,\mu}^\nu(v)=g_{\lambda,\mu}^\nu(v)\frac{a_\nu(v)}{a_\lambda(v)\cdot a_\mu(v)}.
	\end{align*}
	Inspired by \eqref{eq:iterated-Hall}, we can define
	\begin{align}
		\label{eq:Hall-poly}
		f^\lambda_{\lambda_1,\dots,\lambda_n}(v),\quad   g^\lambda_{\lambda_1,\dots,\lambda_n}(v)\in\Z[v,v^{-1}].
	\end{align}

	For $(\alpha,\nu),(\beta,\mu)\in\Z^\I\times\fp$, 
	there exists a polynomial $\boldsymbol{\varphi}^{\lambda,\gamma}_{\mu,\alpha;\nu,\beta}(v)\in\Z[v,v^{-1}]$ such that
	\begin{align}
		\label{eq:multiplication}
		\big([\K_\alpha]\ast[M_q(\mu)]\big)\ast\big([\K_\beta]\ast[M_q(\nu)]\big)=\sum_{\lambda\in\fp,\gamma\in\Z^\I}\boldsymbol{\varphi}^{\lambda,\gamma}_{\mu,\alpha;\nu,\beta}({\sqq})[\K_\gamma]\ast[M_q(\lambda)]
	\end{align}
	in $\widetilde{\ch}(\bfk Q,\btau)$; see \cite[Lemma 9.6]{LW22a}. 
	The generic iHall algebra $\tMHg$ 
	is defined to be the $\Q(v^{1/2})$-linear space with a basis $\{\K_\alpha*\fu_\lambda\mid \alpha\in\Z^\I,\lambda\in\fp\}$ 
	with multiplication
	\begin{align}
		\label{eq:generic-mult}
		(\K_\alpha*\fu_\mu)*(\K_\beta*\fu_\nu)=\sum_{\lambda\in\fp,\gamma\in\Z^\I}\boldsymbol{\varphi}^{\lambda,\gamma}_{\mu,\alpha;\nu,\beta}(v)\K_\gamma*\fu_\lambda.
	\end{align}
	If $\lambda$ is the characteristic function of $\beta\in\Phi^+$, we also denote $\fu_\lambda$ by $\fu_\beta$.

	From Lemmas \ref{lem:Hall-iQG} and \ref{lem:Bridgeland} (see also \cite[Theorem 9.8]{LW22a}), we obtain 
	the isomorphisms of $\Q(v^{1/2})$-algebras
	\begin{alignat*}{2}
		&\widetilde{\psi}:\tUi\longrightarrow \widetilde{\ch}(Q,\btau),\qquad
		&\widetilde{\psi}^{\dbl}:\tU\longrightarrow \widetilde{\ch}(Q^{\rm dbl},\swa).
	\end{alignat*}

	\section{Hall algebra realization of quantum symmetric pairs}
	\label{sec:main}

	In this section, we shall give an embedding of $\widetilde{\ch}(Q,\varrho)$ in $\widetilde{\ch}(Q^{\dbl},\swa)$, which realizes the natural embedding $\imath:\tUi\rightarrow \tU$, and then explain the coideal subalgebra structure of $\tUi$ in $\tU$ in the framework of Hall algebras. We also give the Hopf algebra structure for $\widetilde{\ch}(Q^{\dbl},\swa)$, which is isomorphic to $\tU$ as Hopf algebras.

	\subsection{iHall algebras and iHopf algebras}
	\label{subsec:iHall-iHopf}
	Let $\widetilde{\mathbf{B}}$ be the Borel subalgebra of $\tU$ generated by $E_i$, $K_i$, $i\in\I$. Then $\widetilde{\mathbf{B}}$ is a Hopf subalgebra.
	
	Let $\widetilde{\mathbf{B}}(\bfk Q)$ be the subalgebra of $\widetilde{\ch}(\bfk Q^{\dbl},\swa)$ generated by $[S_i]$, $[\K_{i^\diamond}]$, $i\in\I$; cf. Lemma \ref{lem:Bridgeland}. Its generic version is denoted by $\widetilde{\mathbf{B}}(Q)$. We denote the multiplication of $\widetilde{\mathbf{B}}(\bfk Q)$ (also $\widetilde{\mathbf{B}}(Q)$) by $\cdot$ to avoid confusion. 
	
	For any $\lambda \in \fp$, we denote
	\[ \K_{\lambda}:=\prod_{i \in \I} \K_{i}^{d^\lambda_i} \in \widetilde{\ch}(Q,\varrho), \]
	and for any $\alpha \in \Z^\I$, $\K_{\alpha+\lambda}$ means $\K_{\alpha}\K_{\lambda}$.
	By \cite{Gr95,X97}, we know that  $\widetilde{\mathbf{B}}(Q)$ is a Hopf algebra with the coproduct $\Delta$, the counit $\varepsilon$, and the antipode $S$ given by, for any $\lambda\in\fp$, $\alpha\in\Z^\I$: 
	\begin{align*}
		\Delta(\fu_\lambda)=&\sum_{\mu,\nu\in\fp} v^{\langle \mu,\nu \rangle_Q} f_{\mu,\nu}^\lambda(v) \fu_\mu\cdot \K_{\nu^\diamond}\otimes\fu_\nu,\quad 
		\Delta(\K_{\alpha^\diamond})=\K_{\alpha^\diamond}\otimes \K_{\alpha^\diamond},\quad \forall \alpha\in\Z^\I;
		\\
		\varepsilon(\fu_\lambda)=&\delta_{\lambda,0},\quad\varepsilon(\K_{\alpha^\diamond})=1,\quad S(\K_{\alpha^\diamond})=\K_{\alpha^\diamond}^{-1};
		\\
		S(\fu_\lambda)=& \sum_{m \ge 0}(-1)^m \sum_{\substack{ \mu \in \fp \\ \lambda_1,\dots,\lambda_m \in \fp \backslash \{0\}}} v^{2\sum_{k<l} \EuFo{\lambda_k,\lambda_l}_Q} f^\lambda_{\lambda_1,\dots,\lambda_m}(v) g^\mu_{\lambda_1,\dots,\lambda_m}(v) \K_{\lambda^\diamond}^{-1} \cdot \fu_\mu.
	\end{align*}
	Here and below, for $m=0$ we mean  $f^\lambda_{\lambda_1,\dots,\lambda_m}(v)=g^\lambda_{\lambda_1,\dots,\lambda_m}(v)=\delta_{\lambda,0}$.

	From \cite{Gr95,Rin96}, there is a Hopf pairing $\varphi(-,-)$ on $\widetilde{\mathbf{B}}(Q)$ given by
	\begin{equation}
		\label{eq:hopf-pairing}
		\varphi(\fu_\lambda\cdot \K_{\alpha^{\diamond}},\fu_\mu\cdot \K_{\beta^{\diamond}})=v^{(\alpha,\beta)_Q}\delta_{\lambda,\mu}\cdot a_\lambda(v).
	\end{equation}
	By \cite{CLPRW}, we can give a new associative algebra structure ($\widetilde{\mathbf{B}}(Q)$, $*$) on $\widetilde{\mathbf{B}}(Q)$, with multiplication defined by
	\begin{align}\label{star product}
		a\ast b:=\sum \varphi(\varrho b_{(2)},a_{(1)})\cdot a_{(2)}b_{(1)},\quad \forall a,b\in \widetilde{\mathbf{B}}(Q),
	\end{align}
	where $\Delta(a)=\sum a_{(1)}\otimes a_{(2)}$, $\Delta(b)=\sum b_{(1)}\otimes b_{(2)}$ are the Sweedler notations.
	We denote this new algebra ($\widetilde{\mathbf{B}}(Q)$, $*$) by $\widetilde{\mathbf{B}}^\imath_\varrho(Q)$.
	
	\begin{proposition}
		\label{prop:iHall-ihopf}
		There exists an algebra isomorphism
		\begin{align}
			&\widetilde{\Phi}: \widetilde{\ch}(Q,\varrho)\longrightarrow \widetilde{\mathbf{B}}^\imath_\varrho(Q)  
			\\
			\fu_\lambda*\K_\alpha&\mapsto v^{(\lambda, \alpha)_Q+\EuFo{\alpha,\varrho\alpha}_Q}\fu_{\lambda} \cdot \K_{\varrho\alpha^\diamond},\quad \forall \lambda\in\fp,\alpha\in\Z^\I.
		\end{align}
	\end{proposition}
	
	\begin{proof}
		By Lemma \ref{basis-iHall},  $\widetilde{\Phi}$ is an isomorphism as vector spaces. It remains to prove that $\widetilde{\Phi}$ preserves multiplication. Using Proposition \ref{prop:iHallmult}, in $\widetilde{\ch}(Q,\varrho)$ we have
		\begin{align*}
			(\fu_\lambda * \K_{\alpha}) * (\fu_\mu * \K_{\beta})=&v^{(\varrho \alpha, \mu)_Q-(\alpha,\mu)_Q} \fu_\lambda* \fu_\mu * \K_{\alpha+\beta} 
			\\
			=& \sum_{\nu, \gamma , \delta, \theta \in \fp} v^{\Xi} a_{\delta}(v) f^\lambda_{\delta, \nu}(v) f^\mu_{\gamma, \varrho\delta}(v) g^\theta _{\nu, \gamma}(v)
			\fu_\theta * \K_{\alpha+\beta+\delta},
		\end{align*}
		where
		\[ \Xi= (\varrho \alpha, \mu)_Q-(\alpha,\mu)_Q+\EuFo{\theta,\varrho \delta }_Q-\EuFo{\theta,\delta}_Q-\EuFo{\lambda,\mu}_Q+2\EuFo{\nu, \gamma}_Q.\]
		On the other hand, in $\widetilde{\mathbf{B}}^\imath_\varrho(Q)$
		\begin{align*}
			&(\fu_{\lambda} \cdot \K_{\varrho \alpha^\diamond}) * (\fu_{\mu} \cdot \K_{\varrho\beta^\diamond})
			\\
			&= \sum_{\nu, \gamma , \delta \in \fp} v^{\EuFo{ \delta,  \nu}_Q+ \EuFo{\gamma,\varrho  \delta}_Q+(\beta,\nu+\varrho \alpha)_Q}
			a_{\delta}(v) f^{\lambda}_{\delta, \nu}(v) f^{\mu}_{\gamma, \varrho\delta}(v)
			\fu_{ \nu} \cdot \K_{\varrho\alpha^\diamond} \cdot \fu_{ \gamma} \cdot \K_{\varrho\beta^\diamond+\varrho\delta^\diamond} \\
			&=\sum_{\nu, \gamma , \delta, \theta \in \fp} v^{\Xi'} a_{\delta}(v) f^\lambda_{\delta, \nu}(v) f^\mu_{\gamma, \varrho\delta}(v) g^\theta _{\nu, \gamma}(v) 
			\fu_{\theta} \cdot \K_{\varrho\alpha^\diamond+\varrho\beta^\diamond+\varrho\delta^\diamond},
		\end{align*}
		where
		\[\Xi'=\EuFo{\delta, \nu}_Q+ \EuFo{\gamma, \varrho \delta}_Q+\EuFo{\nu,\gamma}_Q+(\varrho \beta,\varrho \nu+\alpha)_Q+(\alpha,\varrho \gamma)_Q.\]
		The desired result follows by noting that
		\[ \Xi-\Xi'=(\lambda,\alpha)_Q+\EuFo{\alpha,\varrho \alpha}_Q+(\mu,\beta)_Q+\EuFo{\beta,\varrho \beta}_Q-(\theta, \alpha+\beta+\delta)_Q-\EuFo{\alpha+\beta+\delta,\varrho(\alpha+\beta+\delta)}_Q\]
		if $\lambda=\nu+\delta, \; \mu=\gamma+\varrho \delta, \; \theta=\nu+\gamma$.
		The proof is completed.
	\end{proof}
	
	\subsection{Bridgeland's isomorphism as Hopf algebras}
	
	Consider the tensor algebra $\widetilde{\mathbf{B}}( Q)\otimes \widetilde{\mathbf{B}}( Q)$. Following \cite{CLPRW}, we can give a new associative algebra structure on $\widetilde{\mathbf{B}}( Q)\otimes \widetilde{\mathbf{B}}( Q)$, with multiplication given by
	\[(a\otimes b)\ast (c\otimes d)=\sum\varphi (a_{(1)},d_{(2)})\cdot \varphi(c_{(2)},b_{(1)})\cdot a_{(2)}c_{(1)}\otimes b_{(2)}d_{(1)},\qquad\forall a,b,c,d\in \widetilde{\mathbf{B}}( Q).
	\]
	This new algebra is denoted by $D(\widetilde{\mathbf{B}}( Q))$. In fact, $D(\widetilde{\mathbf{B}}( Q))$ is isomorphic to the Drinfeld double of $\widetilde{\mathbf{B}}(Q)$ with respect to $\varphi(-,-)$; see \cite[Proposition 2.8]{CLPRW}. 
	
	\begin{lemma}[{\cite[Proposition 2.7]{CLPRW}}]
		\label{prop:ihopf-Drinfelddouble}
		The   $D(\widetilde{\mathbf{B}}( Q))$ can be made into a Hopf algebra with coproduct, counit and antipode given by 
		\begin{align*}
			\Delta^\imath(a\otimes b)&=\sum\varphi(a_{(2)},b_{(2)})\cdot (a_{(1)}\otimes b_{(3)})\otimes (a_{(3)}\otimes b_{(1)}),
			\\
			\varepsilon^\imath(a\otimes b)&=\varphi(a,S^{-1}(b)),\\
			S^\imath(a\otimes b)&=\sum\varphi(a_{(1)},S^{-1}(b_{(3)}))\cdot\varphi(a_{(2)},b_{(2)})\cdot S(a_{(3)})\otimes S^{-1}(b_{(1)}).
		\end{align*}
	\end{lemma}
	
	For any $N\in\mod(\bfk Q)$, denote by $N^\diamond$ the corresponding module of $\bfk Q^\diamond$. Similarly, we can define $\fu_{\lambda^\diamond}$ for any $\lambda\in\fp$. We denote by $\fu_{\lambda}\oplus\fu_{\mu^\diamond}\in \widetilde{\ch}(Q^{\dbl},\swa)$ the corresponding element of $[M_q(\lambda)\oplus M_q(\mu)^\diamond]\in \widetilde{\ch}(\bfk Q^{\dbl},\swa)$.

	
	\begin{proposition}
		\label{prop: iso of diagonal type}
		There exists an algebra isomorphism
		\begin{align}
			&\widetilde{\Phi}^{\dbl}: \widetilde{\ch}(Q^{\dbl},\swa)\longrightarrow D(\widetilde{\mathbf{B}}( Q))  
			\\
			\fu_{\lambda}\oplus\fu_{\mu^\diamond}&*\K_\alpha*\K_{\beta^\diamond}\mapsto v^{(\mu+\alpha, \lambda+\beta)_Q-(\lambda,\mu)_Q}\fu_\lambda \cdot\K_{\beta^\diamond} \otimes \fu_\mu \cdot \K_{\alpha^\diamond},
		\end{align}
		for any $\lambda,\mu\in\fp,\alpha,\beta\in\Z^\I$. 
	\end{proposition}
	
	\begin{proof}
		By Lemma \ref{basis-iHall}, we only need to show that $\widetilde{\Phi}^{\dbl}$ preserves multiplication. Proposition \ref{prop:iHallmult} gives us the following multiplication formula of $\widetilde{\ch}(Q^{\dbl},\swa)$:
		\begin{align*}
			&(\fu_{\lambda_1} \oplus \fu_{\mu_1^\diamond} * \K_{\alpha_1} * \K_{\beta_1^\diamond}) *  (\fu_{\lambda_2} \oplus \fu_{\mu_2^\diamond} * \K_{\alpha_2} * \K_{\beta_2^\diamond}) \\
			=& \sum_{\lambda, \mu \in \fp} \sum_{\substack{\nu_1, \nu_2, \gamma_1, \gamma_2, \\ \delta_1, \delta_2 \in \fp}} v^{\Xi} a_{\delta_1}(v) a_{\delta_2}(v) f^{\lambda_1}_{\delta_1, \nu_1}(v) f^{\mu_1}_{\delta_2, \nu_2}(v) f^{\lambda_2}_{\gamma_2,\delta_2}(v) f^{\mu_2}_{\gamma_1, \delta_1}(v) g^{\lambda}_{\nu_1, \gamma_2}(v) g^{\mu}_{\nu_2, \gamma_1}(v) 
			\\
			&  \qquad\qquad\qquad\cdot\fu_\lambda \oplus \fu_{\mu^\diamond} * \K_{\alpha_1+ \alpha_2+\delta_1} * \K_{\beta_1^\diamond+ \beta_2^\diamond+\delta_2^\diamond},
		\end{align*}
		where
		\begin{align*}
			\Xi=& \EuFo{ \delta_1, \nu_1 }_Q+\EuFo{\delta_2, \nu_2}_Q+\EuFo{\gamma_1, \delta_1}_Q+\EuFo{\gamma_2,\delta_2}_Q+\EuFo{\nu_1, \gamma_2}_Q+\EuFo{\nu_2,\gamma_1}_Q \\
			&+(\lambda,\mu)_Q-(\lambda+\delta_2, \mu+\delta_1)_Q-(\alpha_1-\beta_1, \lambda_2-\mu_2)_Q.
		\end{align*}
		A direct calculation shows that in $D(\widetilde{\mathbf{B}}(Q))$
		\begin{align*}
			&(\fu_{\lambda_1} \cdot \K_{\beta_1^\diamond} \otimes \fu_{\mu_1} \cdot \K_{\alpha_1^\diamond}) * (\fu_{\lambda_2} \cdot \K_{\beta_2^\diamond}  \otimes \fu_{\mu_2} \cdot \K_{\alpha_2^\diamond}) \\
			=& \sum_{\substack{\nu_1, \nu_2, \gamma_1, \gamma_2, \\ \delta_1, \delta_2 \in \fp,}}
			v^{\EuFo{ \delta_2, \nu_2 }_Q+\EuFo{\gamma_2, \delta_2}_Q+\EuFo{\delta_1, \nu_1}_Q+\EuFo{\gamma_1,\delta_1}_Q+(\nu_1+\beta_1, \alpha_2)_Q+(\nu_2+\alpha_1, \beta_2)_Q}
			a_{\delta_1}(v) a_{\delta_2}(v)
			f^{\mu_1}_{\delta_2, \nu_2}(v) f^{\lambda_2}_{\gamma_2, \delta_2}(v) \\
			&\qquad\qquad\cdot f^{\lambda_1}_{\delta_1, \nu_1}(v) f^{\mu_2}_{\gamma_1, \delta_1}(v) \fu_{\nu_1} \cdot \K_{\beta_1^\diamond} \cdot \fu_{\gamma_2} \cdot \K_{\beta_2^\diamond+\delta_2^\diamond} 
			\otimes \fu_{\nu_2} \cdot \K_{\alpha_1^\diamond} \cdot \fu_{\gamma_1} \cdot \K_{\alpha_2^\diamond+\delta_1^\diamond} \\
			=& \sum_{\lambda, \mu \in \fp} \sum_{\substack{\nu_1, \nu_2, \gamma_1, \gamma_2, \\ \delta_1, \delta_2 \in \fp}} v^{\Xi'} a_{\delta_1}(v) a_{\delta_2}(v)
			f^{\mu_1}_{\delta_2, \nu_2}(v) f^{\lambda_2}_{\gamma_2, \delta_2}(v) f^{\lambda_1}_{\delta_1, \nu_1}(v) f^{\mu_2}_{\gamma_1, \delta_1}(v) g^{\lambda}_{\nu_1, \gamma_2}(v) g^{\mu}_{\nu_2, \gamma_1}(v)
			\\
			&\qquad\qquad\qquad\cdot\fu_{\lambda} \cdot \K_{\beta_1^\diamond+\beta_2^\diamond+\delta_2^\diamond} \otimes \fu_{\mu} \cdot \K_{\alpha_1^\diamond+\alpha_2^\diamond+\delta_1^\diamond}, 
		\end{align*}
		where
		\begin{align*}
			\Xi'=&\EuFo{ \delta_1, \nu_1 }_Q+\EuFo{\delta_2, \nu_2}_Q+\EuFo{\gamma_1, \delta_1}_Q+\EuFo{\gamma_2,\delta_2}_Q+\EuFo{\nu_1, \gamma_2}_Q+\EuFo{\nu_2,\gamma_1}_Q \\
			&+(\alpha_1, \gamma_1)_Q+(\beta_1,\gamma_2)_Q +(\nu_2+\alpha_1, \beta_2)_Q+(\nu_1+\beta_1, \alpha_2)_Q. 
		\end{align*}
		Since
		\begin{align*}
			\Xi-\Xi'=&(\mu_1+\alpha_1, \lambda_1+\beta_1)_Q-(\lambda_1,\mu_1)_Q+(\mu_2+\alpha_2, \lambda_2+\beta_2)_Q-(\lambda_2,\mu_2)_Q \\
			&-(\mu+\alpha_1+\alpha_2+\delta_1, \lambda+\beta_1+\beta_2+\delta_2)_Q+(\lambda,\mu)_Q
		\end{align*}
		if
		\[ \lambda_1=\delta_1+\nu_1, \; \mu_1=\delta_2+\nu_2, \; \lambda_2=\gamma_2+\delta_2, \; \mu_2=\gamma_1+\delta_1, \; \lambda=\nu_1+\gamma_2, \; \mu=\nu_2+\gamma_1, \]
		the proof is completed.
	\end{proof}

	Now, we can give the Hopf algebra structure of $\widetilde{\ch}( Q^{\dbl},\swa)$.
	
	\begin{theorem}
		\label{thm:Hopf-iHall}
		The $\widetilde{\ch}(Q^{\dbl},\swa)$ is a Hopf algebra with the coproduct, the counit and the antipode given by, for any $\lambda,\mu\in\fp$, $\alpha,\beta\in\Z^\I$,
		\begin{align}
			\Delta(\K_\alpha)=&\K_\alpha\otimes \K_\alpha,\qquad \Delta(\K_{\alpha^\diamond})=\K_{\alpha^\diamond}\otimes \K_{\alpha^\diamond},
			\\
			\label{coprod}
			\Delta(\fu_\lambda\oplus \fu_{\mu^\diamond})=& \sum_{\substack{\lambda_1,\lambda_2, \\ \mu_1,\mu_2,\nu \in \fp}} v^{\EuFo{\lambda_1, \lambda_2+\nu}_Q+\EuFo{\mu_1,\mu_2+\nu}_Q-\EuFo{\lambda_2+\mu_2,\nu}_Q-(\lambda_2,\mu_2)_Q} \cdot a_\nu(v) \cdot f^\lambda_{\lambda_1,\nu,\lambda_2} (v)\\ \notag
			&\qquad\;\; \quad \cdot f^\mu_{\mu_1,\nu,\mu_2}(v)   (\fu_{\lambda_1} \oplus \fu_{\mu_2^\diamond}*\K_{\lambda_2^\diamond+\nu^\diamond}) \otimes(\fu_{\lambda_2} \oplus \fu_{\mu_1^\diamond} *\K_{\mu_2+\nu}),
			\\
			\varepsilon(\K_\alpha)=&1=\varepsilon(\K_{\alpha^\diamond}),
			\\
			\label{varepsilon}
			\varepsilon(\fu_\lambda\oplus \fu_{\mu^\diamond})=&a_\lambda(v)\cdot\sum_{m \ge 0} (-1)^m \sum_{\mu_1\dots,\mu_m \in \fp \backslash \{0\}} v^{\sum_{k<l}(\mu_k,\mu_l)_Q} g^\lambda_{\mu_m,\dots,\mu_1}(v)\cdot f^\mu_{\mu_1,\dots,\mu_m}(v),
			\\
			S(\K_\alpha)=&\K_\alpha^{-1},\qquad  S(\K_{\alpha^\diamond})=\K_{\alpha^\diamond}^{-1},
			\\
			\label{antipode}
			S(\fu_\lambda\oplus \fu_{\mu^\diamond})=& \sum_{\nu, \omega,\pi \in \fp}
			\sum_{r, s, t \ge 0} \sum_{\gamma_1,\dots,\gamma_r \in \fp \backslash \{ 0 \}}  \sum_{\delta_1,\dots,\delta_s \in \fp \backslash \{ 0 \}}  \sum_{\theta_1,\dots,\theta_t \in \fp \backslash \{ 0 \}}  (-1)^{r+s+t}v^{\EuFo{\lambda-\omega,\omega-\pi}_Q}\\ \notag
			& \cdot v^{-(\omega,\omega-\pi)_Q+\sum_{k<l}(\gamma_k,\gamma_l)_Q+2\sum_{k<l}\EuFo{\delta_k,\delta_l}_Q+\sum_{k<l}(\theta_k,\theta_l)_Q}\cdot a_\nu(v) \cdot a_{\gamma_1}(v) \cdots a_{\gamma_r}(v)\\\notag
			& \cdot   f^\lambda_{\gamma_r,\dots,\gamma_1,\nu,\delta_1,\dots,\delta_s}(v) f^\mu_{\theta_1,\dots,\theta_t,\nu,\gamma_1,\dots,\gamma_r}(v) g^\omega_{\delta_1,\dots,\delta_s}(v) g^\pi_{\theta_t,\dots,\theta_1}(v) \fu_\omega \oplus \fu_{\pi^\diamond}*\K_{\mu}^{-1} *\K_{\omega^\diamond}^{-1} .
		\end{align}
	\end{theorem}
	
	\begin{proof}
		One only need to verify that the Hopf structure of $\widetilde{\ch}(Q^{\dbl},\swa)$ defined here is deduced from $D(\widetilde{\mathbf{B}}(Q))$ via the algebra isomorphism $\widetilde{\Phi}^{\dbl}$ in Proposition \ref{prop: iso of diagonal type}.
		
		For \eqref{coprod} we compute
		\begin{align*}
			\Delta^\imath \circ\widetilde{\Phi}^{\dbl} (\fu_\lambda \oplus \fu_{\mu^\diamond})=& \Delta^\imath(\fu_{\lambda} \otimes \fu_\mu) \\
			=& \sum_{\substack{\lambda_1,\lambda_2, \\\mu_1,\mu_2,\nu  \in \fp}} v^{\EuFo{\lambda_1,\lambda_2+\nu}_Q+\EuFo{\nu,\lambda_2}_Q+\EuFo{\mu_1,\mu_2+\nu}_Q+\EuFo{\nu,\mu_2}_Q+(\mu_2,\lambda_2)_Q}\cdot a_\nu(v) \\
			&\cdot f^\lambda_{\lambda_1,\nu,\lambda_2}(v)\cdot f^\mu_{\mu_1,\nu,\mu_2}(v) (\fu_{\lambda_1} \cdot \K_{\lambda_2^\diamond+\nu^\diamond} \otimes \fu_{\mu_2}) \otimes (\fu_{\lambda_2} \otimes \fu_{\mu_1} \cdot \K_{\mu_2^\diamond+\nu^\diamond}).
		\end{align*}
		Then \eqref{coprod} follows by applying
		$(\widetilde{\Phi}^{\dbl}\otimes \widetilde{\Phi}^{\dbl})^{-1}$ to the above formula.
		
		The inverse of antipode $S$ of $\widetilde{\mathbf{B}}(Q)$ is as follows (see \cite{X97}):
		\[ S^{-1}(\fu_\lambda)=\sum_{m \ge 0} (-1)^m \sum_{\substack{\pi \in \fp \\ \lambda_1,\dots,\lambda_m \in \fp \backslash \{0\}}} v^{\sum_{k<l}(\lambda_k,\lambda_l)_Q} f^\lambda_{\lambda_1,\dots,\lambda_m}(v) g^\pi_{\lambda_m,\dots,\lambda_1}(v) \fu_{\pi} \cdot \K^{-1}_{\lambda^\diamond}. \]
		For \eqref{varepsilon} we have
		\begin{align*}
			\varepsilon^\imath\circ \widetilde{\Phi}^{\dbl}(\fu_\lambda \oplus \fu_{\mu^\diamond})=&\varepsilon^\imath(\fu_{\lambda} \otimes \fu_\mu) \\
			=&\sum_{m \ge 0} (-1)^m \sum_{\substack{\pi \in \fp \\ \mu_1,\dots,\mu_m \in \fp \backslash \{0\}}} v^{\sum_{k<l}(\mu_k,\mu_l)_Q} f^\mu_{\mu_1,\dots,\mu_m}(v) g^\pi_{\mu_m,\dots,\mu_1}(v)  \varphi(\fu_{\lambda},\fu_{\pi} \cdot \K^{-1}_{\mu^\diamond}) \\
			=& a_\lambda(v) \sum_{m \ge 0} (-1)^m \sum_{\mu_1,\dots,\mu_m \in \fp \backslash \{0\}} v^{\sum_{k<l}(\mu_k,\mu_l)_Q}g^\lambda_{\mu_m,\dots,\mu_1}(v) f^\mu_{\mu_1,\dots,\mu_m}(v).
		\end{align*}
		Then \eqref{varepsilon} follows since
		$\varepsilon(\fu_\lambda \oplus \fu_{\mu^\diamond})=\varepsilon^\imath\circ \widetilde{\Phi}^{\dbl}(\fu_\lambda \oplus \fu_{\mu^\diamond})$.

		Finally, for \eqref{antipode} a calculation shows that
		\begin{align*}
			&S^\imath \circ \widetilde{\Phi}^{\dbl}(\fu_\lambda \oplus \fu_{\mu^\diamond})= S^\imath(\fu_{\lambda} \otimes \fu_\mu) \\
			=&  \sum_{\substack{\lambda_1,\lambda_2, \\\mu_1,\mu_2,\nu  \in \fp}}v^{\EuFo{\lambda_1,\lambda_2+\nu}_Q+\EuFo{\nu,\lambda_2}_Q+\EuFo{\mu_1,\mu_2+\nu}_Q+\EuFo{\nu,\mu_2}_Q+(\mu_2,\lambda_2)_Q} \\
			& \cdot a_\nu(v) f^\lambda_{\lambda_1,\nu,\lambda_2}(v) f^\mu_{\mu_1,\nu,\mu_2}(v) \varphi(\fu_{\lambda_1} \cdot \K_{\lambda_2^\diamond+\nu^\diamond},S^{-1}(\fu_{\mu_2})) \cdot S(\fu_{\lambda_2}) \otimes S^{-1}(\fu_{\mu_1} \cdot \K_{\mu_2^\diamond+\nu^\diamond})\\
			=& \sum_{\substack{\lambda_1,\lambda_2, \mu_1,\mu_2, \\ \nu, \omega,\pi \in \fp}}
			\sum_{r, s, t \ge 0} \sum_{\gamma_1,\dots,\gamma_r \in \fp \backslash \{ 0 \}}  \sum_{\delta_1,\dots,\delta_s \in \fp \backslash \{ 0 \}}  \sum_{\theta_1,\dots,\theta_t \in \fp \backslash \{ 0 \}}  (-1)^{r+s+t} v^{(\mu_2,\lambda_2)_Q} \\
			&\cdot v^{\EuFo{\lambda_1,\lambda_2+\nu}_Q+\EuFo{\nu,\lambda_2}_Q+\EuFo{\mu_1,\mu_2+\nu}_Q+\EuFo{\nu,\mu_2}_Q+\sum_{k<l}(\gamma_k,\gamma_l)_Q+2\sum_{k<l}\EuFo{\delta_k,\delta_l}_Q+\sum_{k<l}(\theta_k,\theta_l)_Q} \\
			&\cdot a_\nu(v)\cdot f^\lambda_{\lambda_1,\nu,\lambda_2}(v) \cdot f^\mu_{\mu_1,\nu,\mu_2}(v) v^{-(\lambda_2+\nu,\mu_2)_Q} a_{\lambda_1}(v) \cdot f^{\mu_2}_{\gamma_1,\dots,\gamma_r}(v) \cdot g^{\lambda_1}_{\gamma_r,\dots,\gamma_1}(v) \cdot f^{\lambda_2}_{\delta_1,\dots,\delta_s}(v) 
			\\
			&\cdot g^\omega_{\delta_1,\dots,\delta_s}(v)\cdot f^{\mu_1}_{\theta_1,\dots,\theta_t}(v)\cdot g^\pi_{\theta_t,\dots,\theta_1}(v)\cdot
			\K_{\lambda_2^\diamond}^{-1} \cdot \fu_{\omega} \otimes \K_{\mu_2^\diamond+\nu^\diamond}^{-1} \cdot \fu_\pi \cdot \K_{\mu_1^\diamond}^{-1} \\
			=& \sum_{\nu, \omega,\pi \in \fp}
			\sum_{r, s, t \ge 0} \sum_{\gamma_1,\dots,\gamma_r \in \fp \backslash \{ 0 \}}  \sum_{\delta_1,\dots,\delta_s \in \fp \backslash \{ 0 \}}  \sum_{\theta_1,\dots,\theta_t \in \fp \backslash \{ 0 \}}  (-1)^{r+s+t}\\ \notag
			& \cdot v^{\EuFo{\lambda-\omega,\omega-\pi}_Q-(\omega,\omega)_Q+\sum_{k<l}(\gamma_k,\gamma_l)_Q+2\sum_{k<l}\EuFo{\delta_k,\delta_l}_Q+\sum_{k<l}(\theta_k,\theta_l)_Q} a_\nu(v)  \cdot a_{\gamma_1}(v) \cdots a_{\gamma_r}(v)\\\notag
			& \cdot f^\lambda_{\gamma_r,\dots,\gamma_1,\nu,\delta_1,\dots,\delta_s}(v)\cdot f^\mu_{\theta_1,\dots,\theta_t,\nu,\gamma_1,\dots,\gamma_r}(v)\cdot    g^\omega_{\delta_1,\dots,\delta_s} (v)\cdot g^\pi_{\theta_t,\dots,\theta_1}(v) \cdot \fu_{\omega} \cdot \K_{\omega^\diamond}^{-1} \otimes \fu_\pi \cdot \K^{-1}_{\mu^\diamond}.
		\end{align*}
		Then \eqref{antipode} follows by applying $(\widetilde{\Phi}^{\dbl})^{-1}$ to the above formula.
	\end{proof}
	
	\begin{remark}
		The Hopf algebra structure in Theorem \ref{thm:Hopf-iHall} holds for arbitrary quivers $Q$ (without the assumption of Dynkin quivers). 
		It is remarkable that \eqref{coprod} has been obtained in \cite[Theorem 3.4]{LiP25} by a long direct computation. Our method is much simpler; cf. \cite[Proposition 2.7]{CLPRW}.
	\end{remark}

	\begin{corollary}
		\label{cor:Bridgeland-Hopf}
		The algebra isomorphism $\widetilde{\psi}^{\dbl}:\tU\stackrel{\simeq}{\longrightarrow} \widetilde{\ch}(Q^{\rm dbl},\swa)$ given in Lemma \ref{lem:Bridgeland} is an isomorphism as Hopf algebras.
	\end{corollary}
	
	\begin{proof}
		To avoid confusion, we denote the comultiplication and the counit of $\widetilde{\ch}(Q^{\rm dbl},\swa)$ by $\Delta_{\widetilde{\ch}}$ and $\varepsilon_{\widetilde{\ch}}$ respectively. Since
		\[\varepsilon_{\widetilde{\ch}}(\fu_{\lambda})=\varepsilon_{\widetilde{\ch}}(\fu_{\lambda^\diamond})=0, \quad  \varepsilon_{\widetilde{\ch}}(\K_{\alpha})=\varepsilon_{\widetilde{\ch}}(\K_{\alpha^\diamond})=1\] for any $\lambda \in \fp \backslash\{0\}$ and $\alpha \in \Z^\I$, $\widetilde{\psi}^{\dbl}$ preserves counit. It remains to prove $\widetilde{\psi}^{\dbl}$ preserves comultiplication. Indeed, for any $i \in \I$ we have
		\[ (\widetilde{\psi}^{\dbl} \otimes \widetilde{\psi}^{\dbl}) \circ \Delta_{\tU}(E_i)=v^{-\frac{1}{2}} (\fu_{\alpha_i} \otimes 1+ \K_{\alpha_i^\diamond}\otimes \fu_{\alpha_i}), \]
		\[ \Delta_{\widetilde{\ch}} \circ \widetilde{\psi}^{\dbl}(E_i)=v^{-\frac{1}{2}}(\fu_{\alpha_i} \otimes 1+ \K_{\alpha_i^\diamond} \otimes \fu_{\alpha_i}). \]
		Here $\Delta_{\tU}$ is the coproduct of $\tU$. 
		Thus, $ (\widetilde{\psi}^{\dbl} \otimes \widetilde{\psi}^{\dbl}) \circ \Delta_{\widetilde{\ch}}(E_i)=\Delta_{\widetilde{\ch}} \circ \widetilde{\psi}^{\dbl}(E_i)$. A similar argument shows that $ (\widetilde{\psi}^{\dbl} \otimes \widetilde{\psi}^{\dbl}) \circ \Delta_{\tU}$ and $\Delta_{\widetilde{\ch}} \circ \widetilde{\psi}^{\dbl}(E_i)$ also agree on $F_i, K_i$ and $K_i'$. In other words, $\widetilde{\psi}^{\dbl}$ preserves the comultiplication of all generators of $\tU$. Therefore $\widetilde{\psi}^{\dbl}$ is a bialgebra isomorphism, and the desired result follows from the fact that any bialgebra homomorphism of Hopf algebras is a Hopf homomorphism.
	\end{proof}

	\subsection{Quantum symmetric pairs via Hall algebras}
	
	
	For the set $\Phi^+$ of positive root, by considering $\varrho$, we denote by \begin{align}
		\Phi^+_0:=\{\beta\in\Phi^+\mid\varrho(\beta)=\beta\},\qquad \Phi^+_1:=\{\beta\in\Phi^+\mid\varrho(\beta)\neq\beta\}.
	\end{align}

	For any $M_q(\lambda)\in\mod(\bfk Q)$, we know
	$$M_q(\lambda)\cong \bigoplus_{\beta\in\Phi^+} M_q(\beta)^{\oplus \lambda(\beta)}.$$
	We define a subset $\cf$ of $\fp$ to be
	\begin{align}
		\cf:=\{\lambda\in\fp\mid  \varrho(\lambda)=\lambda,2\mid \lambda(\beta), \forall \beta\in\Phi^+_0\}.
	\end{align}
	For any $\lambda\in\cf$, 
	we define $\lambda_f,\lambda_c\in\fp$ by
	$$\lambda_f(\beta)=\begin{cases}
		\frac{\lambda(\beta)}{2}&\text{ if }\beta\in\Phi^+_0,\\
		0&\text{ otherwise;}
	\end{cases}
	\qquad 
	\lambda_c(\beta)= \begin{cases}
		\lambda(\beta)&\text{ if }\beta\in\Phi^+_1,\\
		0&\text{ otherwise.}
	\end{cases}$$
	By definition, we have 
	$M_q(\lambda)\cong M_q(\lambda_f)\oplus M_q(\lambda_f)\oplus M_q(\lambda_c)$ for any $\lambda\in\cf$. In particular, we have
	$\varrho(M_q(\lambda_c))=M_q(\lambda_c)$, and any indecomposable direct summand of $M_q(\lambda_c)$ is not fixed by $\varrho$.
	
	\begin{lemma}
		\label{lem:sqrt-aut}
		For any $\lambda\in\cf$, there exists $b_{\lambda_c}(v)\in\Z[v,v^{-1}]$ such that $a_{\lambda_c}(v)=b_{\lambda_c}(v)^2$.
	\end{lemma}
	
	Before proceeding with the proof, let us recall some results on $Q$-admissible sequences; see \cite{DDPW,LW21a}. Recall from \S\ref{sub:generic} that $M_q(\beta)$ is the indecomposable $\bfk Q$-module such that $\dimv M_q(\beta)=\beta$ for each $\beta\in\Phi^+$. An ordering $\beta_1,\dots,\beta_N$ of the positive roots in $\Phi^+$ is called to be {\em $Q$-admissible} if $\Ext^1_{\bfk Q}(M_q(\beta_r),M_q(\beta_t))\neq 0$ implies $r>t$ (or $\Hom_{\bfk Q}(M_q(\beta_r),M_q(\beta_t))\neq0$ only if $r\leq t$). Note that this definition does not depend on the base field $\bfk=\F_q$. A $Q$-admissible sequence exists, see, e.g., \cite[Corollary 3.34]{DDPW}.

	From \cite[Theorem 8.3]{LW21a}, we can get a $Q$-admissible sequence of all the roots in $\Phi^+$:  
	\begin{align}
		\label{eq:Qadmissible}
		\beta_{1}, \btau(\beta_1), \beta_{2}, \btau(\beta_2), \dots,\beta_{N_\imath },\btau(\beta_{N_\imath}).
	\end{align}
	(By convention the redundant $\btau(\beta_j)$ is omitted here and below whenever $\btau(\beta_j)=\beta_j$.) In particular, we have
	$$\Hom_{\bfk Q}(M_q(\beta), M_q(\varrho \beta))=\begin{cases}
		\bfk & \text{ if }\btau(\beta)= \beta,\\
		0&\text{otherwise}.
	\end{cases}$$

	For any $\beta\in\Phi^+$, we have 
	$$|\Aut(M_q(m\beta))|=(q^m-1)(q^m-q)\cdots(q^m-q^{m-1}),$$
	and then
	\begin{align}
		\label{function a}
		a_{m \beta}(v)=v^{\frac{m(3m-1)}{2}} (v-v^{-1})^m [m]!.
	\end{align}
	
	\begin{proof}[Proof of Lemma \ref{lem:sqrt-aut}]
		Without loss of generality, 
		we suppose $\lambda=\lambda_c$. We also assume that $ \beta_{1}, \btau(\beta_1), \beta_{2}, \btau(\beta_2), \dots,\beta_{m },\btau(\beta_{m})$ is a $Q$-admissible sequence of the positive roots in $\Phi^+_1$; cf. \eqref{eq:Qadmissible}. 
		Write
		\[ M_q(\lambda)=\bigoplus_{k=1}^m (M_q(\beta_k) \oplus M_q(\varrho \beta_k))^{\oplus r_k}. \]
		Since  \begin{align*}&\Hom_{\bfk Q}(M_q(\beta_k), M_q(\beta_l))=\Hom_{\bfk Q}(M_q(\varrho \beta_k), M_q(\beta_l))\\
			&=\Hom_{\bfk Q}(M_q(\beta_k), M_q(\varrho \beta_l))=\Hom_{\bfk Q}(M_q(\varrho \beta_k), M_q(\varrho \beta_l))=0
		\end{align*}
		for any $k >l$, any endomorphism of $M_q(\lambda)$ can be represented by a block lower triangular matrix, whose $(l,k)$-block is a homomorphism from $(M_q(\beta_k) \oplus M_q(\varrho \beta_k))^{\oplus r_k}$ to $(M_q(\beta_l) \oplus M_q(\varrho \beta_l))^{\oplus r_l}$. Therefore, an endomorphism $\phi$ of $M_q(\lambda)$ is an automorphism if and only if each $(k,k)$-block of $\phi$ is an automorphism of $(M_q(\beta_k) \oplus M_q(\varrho \beta_k))^{\oplus r_k}$. 
		
		Recall $\Hom_{\bfk Q}(M_q(\beta), M_q(\varrho \beta))= 0$ if $\beta \in \Phi^+_1$. So
		\[ a_{M_q(\lambda)}=\prod_{k < l} q^{\dim \Hom_{\bfk Q}((M_q(\beta_k) \oplus M_q(\varrho \beta_k))^{\oplus r_k}, (M_q(\beta_l) \oplus M_q(\varrho \beta_l))^{\oplus r_l})}
		\prod_{k=1}^m a_{M_q(r_k \beta_k)} a_{M_q(r_k \varrho \beta_k)}.\]
		Note that $a_{M_q(r_k \beta_k)}=a_{M_q(r_k \varrho \beta_k)}$, we set
		\[ b_{\lambda}(v)=\prod_{k<l} v^{2r_kr_l \dim_\bfk \Hom(M_q({\beta_k}), M_q({\beta_l}) \oplus M_q({\varrho \beta_l}))} \prod_{k=1}^{m} a_{r_k \beta_k}(v). \]
		Then $b_{\lambda_c}(v)^2=a_{\lambda_c}(v)=a_\lambda(v)$ since $\F_q$ is arbitrary. 
		By using (\ref{function a}), we have $b_{\lambda_c}(v) \in \Z[v,v^{-1}]$ (even in $\Z[v]$). 
	\end{proof}

	\begin{definition}[\cite{CLPRW}]
		A linear map $\chi:\widetilde{\mathbf{B}}( Q)\longrightarrow \Q(v^{\frac{1}{2}})$ is called $\varrho$-twisted compatible if 
		$\chi(a \cdot b)=\sum\chi(a_{(1)})\chi(b_{(2)})\varphi(\varrho (a_{(2)}),b_{(1)})$ holds for all $a,b\in \widetilde{\mathbf{B}}( Q)$. 
	\end{definition}
	
	Let us give some properties of  a $\varrho$-twisted compatible linear map $\chi:\widetilde{\mathbf{B}}( Q)\longrightarrow \Q(v^{\frac{1}{2}})$. 
	\begin{lemma}
		\label{lem: properties of chi}
		Let $\chi:\widetilde{\mathbf{B}}( Q)\longrightarrow \Q(v^{\frac{1}{2}})$ be a $\varrho$-twisted compatible linear map. Then we have the following.
		\begin{enumerate}
			\item 
			For any $\lambda,\mu\in\fp$, if $\Hom(M_q(\mu),\varrho M_q(\lambda))=0$, then 
			\begin{align*}
				\chi(\fu_\lambda\cdot\fu_\mu)=\chi(\fu_\lambda)\cdot \chi(\fu_\mu).
			\end{align*}
			\item For any $\beta\in\Phi_0^+$, we have, for $m\geq2$,
			\begin{align*}
				\chi(\fu_{m \beta})=v^{3m-4}(v-v^{-1})[m-1] \cdot\chi(1)\cdot \chi(\fu_{(m-2) \beta} \cdot \K_{\beta^\diamond})+v^{m-1}\chi(\fu_{(m-1)\beta}) \cdot\chi(\fu_\beta).
			\end{align*}
			\item For any $\beta\in\Phi_1^+$, we have
			\begin{align*}
				\chi(\fu_{m_1 \beta+m_2 \varrho \beta})=&v^{2m_1+m_2-2}(v-v^{-1})[m_1]\cdot \chi(1) \cdot \chi(\fu_{(m_1-1)\beta+(m_2-1)  \varrho \beta} \cdot \K_{\beta^\diamond})\\
				&+v^{m_2-1}\chi(\fu_{m_1 \beta+(m_2-1) \varrho \beta}) \cdot\chi(\fu_{\varrho \beta}),\quad \forall m_1,m_2 \ge 1.
			\end{align*}
			\item Assume $\chi(\fu_{\alpha_i})=0$ for all simple roots $\alpha_i$. For any $h=\fu_{\alpha_{i_1}}\cdots \fu_{\alpha_{i_m}}$, if $\wt(h)=\alpha_{i_1}+\cdots +\alpha_{i_m}\notin \sum_{i \in \I} \Z(\alpha_i +\alpha_{\varrho i})$, then $\chi(h)=0$. In particular, we have $\chi(\fu_\beta)=0$ for any $\beta\in\Phi^+$ in this case. 
		\end{enumerate}
	\end{lemma}

	\begin{proof}
		(1) We have in $\widetilde{\mathbf{B}}(Q)$:
		\begin{align*}
			\Delta(\fu_\lambda)=&\sum_{\lambda',\lambda''\in\fp} v^{\langle \lambda',\lambda'' \rangle_Q} f_{\lambda',\lambda''}^\lambda(v) \fu_{\lambda'}\cdot \K_{\lambda''^\diamond}\otimes\fu_{\lambda''},
			\\
			\Delta(\fu_\mu)=&\sum_{\mu',\mu''\in\fp} v^{\langle \mu',\mu'' \rangle_Q} f_{\mu',\mu''}^\mu(v) \fu_{\mu'}\cdot \K_{\mu''^\diamond}\otimes\fu_{\mu''}.
		\end{align*}
		Note that $f_{\lambda',\lambda''}^\lambda(v)\neq0$ only if $M_q(\lambda'')\subseteq M_q(\lambda)$, and $f_{\mu',\mu''}^\mu(v)\neq0$ only if $M_q(\mu')$ is a quotient of $M_q(\mu)$.
		If $\varphi(\varrho \fu_{\lambda''},\fu_{\mu'}\cdot \K_{\mu''^\diamond})\neq 0$, then $\varrho M_q(\lambda'')\cong M_q(\mu')$ by \eqref{eq:hopf-pairing}. Then we have a morphism $M_q(\mu)\twoheadrightarrow M_q(\mu')\stackrel{\cong}{\longrightarrow} \varrho M_q(\lambda'')\hookrightarrow \varrho M_q(\lambda)$. By assumption, we have  $M_q(\lambda'')=0=M_q(\mu')$. So we have
		\begin{align*}
			\chi(\fu_\lambda\cdot \fu_\mu)=\chi(\fu_\lambda)\cdot \chi(\fu_\mu).
		\end{align*}
		
		(2) Since $\Ext^1(M_q(\beta),M_q(\beta))=0$ we have $\fu_{m \beta}=v^{\frac{m(m-1)}{2}} \fu_\beta^m$. Then
		\[ \chi(\fu_{m \beta})=v^{m-1} \chi(\fu_{(m-1) \beta} \cdot \fu_\beta).\]
		In $\widetilde{\mathbf{B}}(Q)$, by definition, 
		\begin{align*}
			\Delta(\fu_{(m-1) \beta})=&\sum_{\lambda',\lambda''\in\fp} v^{\langle \lambda',\lambda'' \rangle_Q} f_{\lambda',\lambda''}^{(m-1) \beta}(v) \fu_{\lambda'}\cdot \K_{\lambda''^\diamond}\otimes\fu_{\lambda''},
			\\
			\Delta(\fu_\beta)=&\sum_{\mu',\mu''\in\fp} v^{\langle \mu',\mu'' \rangle_Q} f_{\mu',\mu''}^\beta(v) \fu_{\mu'}\cdot \K_{\mu''^\diamond}\otimes\fu_{\mu''}.
		\end{align*}
		If $f_{\lambda',\lambda''}^{(m-1) \beta}(v), f_{\mu',\mu''}^\beta(v), \varphi(\varrho \fu_{\lambda''}, \fu_{\mu'} \cdot \K_{\mu''^\diamond})$ are all nonzero and $\mu' \neq 0$, then there exist nonzero homomorphisms $M_q(\beta) \twoheadrightarrow M_q(\mu')$ and $M_q(\mu') \rightarrow M_q(\beta)$. So in this case $\lambda''=\mu'=\beta$, and moreover $\lambda'=(m-2)\beta, \mu''=0$. Therefore,
		\begin{align*}
			\chi(\fu_{m \beta})=&v^{m-1+(m-2)\EuFo{\beta,\beta}_Q+m-2} [m-1] \chi(\fu_{(m-2) \beta} \cdot \K_{\beta^\diamond}) \chi(1) \varphi(\fu_\beta, \fu_\beta)+v^{m-1}\chi(\fu_{(m-1)\beta}) \chi(\fu_\beta) \\
			=& v^{3m-4}(v-v^{-1})[m-1] \chi(1) \chi(\fu_{(m-2) \beta} \cdot \K_{\beta^\diamond})+v^{m-1}\chi(\fu_{(m-1)\beta}) \chi(\fu_\beta).
		\end{align*}
		
		(3) An argument similar to (2) gives the desired result.
		
		(4) We prove this by induction on $m$. If $m=1$, there is nothing to prove. For $m>1$, recall that $\Delta(\fu_{\alpha_i})=\fu_{\alpha_i} \otimes 1+\K_{\alpha_i^\diamond} \otimes \fu_{\alpha_i}$ for any simple root $\alpha_i$. Since $\chi$ is $\varrho$-twisted compatible, we have
		\begin{align*}
			\chi(h)=&\chi(\fu_{\alpha_{i_1}} \cdots \fu_{\alpha_{i_{m-1}}})\cdot \chi(\fu_{\alpha_{i_m}})\cdot\varphi(1, \K_{\alpha_{i_m}^\diamond})\\
			&+\sum_{k=1}^{m-1}v^{(\alpha_{i_k}, \alpha_{i_{k+1}}+\cdots+\alpha_{i_{m-1}})_Q} \chi(\fu_{\alpha_{i_1}} \cdots \fu_{\alpha_{i_{k-1}}} \fu_{\alpha_{i_{k+1}}}\cdots \fu_{\alpha_{i_{m-1}}} \cdot \K_{\alpha_{i_k}^\diamond}) \chi(1) \cdot\varphi(\fu_{\varrho \alpha_{i_k}}, \fu_{\alpha_{i_m}}) \\
			=& \chi(1) (v^2-1) \sum_{k=1}^{m-1} \delta_{\varrho \alpha_{i_k},\alpha_{i_m}}v^{(\alpha_{i_k}, \alpha_{i_{k+1}}+\cdots+\alpha_{i_{m-1}})_Q}  \chi(h'_k) \cdot \chi(\K_{\alpha_{i_k}^\diamond}),
		\end{align*}
		where $h_k'= \fu_{\alpha_{i_1}} \cdots \fu_{\alpha_{i_{k-1}}} \fu_{\alpha_{i_{k+1}}}\cdots \fu_{\alpha_{i_{m-1}}}$. It is clear that if $\varrho \alpha_{i_k}=\alpha_{i_m}$ then $\textup{wt}(h'_k)=\textup{wt}(h)-(\alpha_{i_k}+\varrho \alpha_{i_k}) \notin \sum_{i \in \I}\Z(\alpha_i + \alpha_{\varrho i})$. The proof is completed by using the inductive hypothesis.
		%
		
		For the last statement, for any $\beta\in\Phi^+$,  if $\varrho=\Id$, we can see $\beta\notin \sum_{i\in\I}2\N\alpha_i$. For $\varrho\neq \Id$, the sole possibility of $\beta\in \sum_{i \in \I} \Z(\alpha_i +\alpha_{\varrho i})$ arises in type $A_{2r}$, a case excluded from our setting of iquivers;   
		cf. \eqref{diag: A}--\eqref{diag: E}. 
	\end{proof}
	
	Let $\I_0$ be the subset of $\I$ formed by $i\in\I$ such that $\varrho(i)=i$. 
	We define
	\begin{align}
		\{\lambda,\mu\}:=\dim_\bfk\Hom(M_q(\lambda),M_q(\mu))+\dim_\bfk\Hom(M_q(\mu),M_q(\lambda)),\quad \forall \lambda,\mu\in\fp.
	\end{align}
	
	Let us define a linear map 
	$$\chi: \widetilde{\mathbf{B}}( Q)\longrightarrow \Q(v^{\frac{1}{2}})$$
	to be such that 
	\begin{align}
		\label{eq:chi1}
		\chi(\K_{\beta^\diamond})=& v^{\langle\beta,\varrho(\beta)\rangle_Q+\sum_{i\in\I_0} d_i^\beta},
		\\
		\label{eq:chi2}
		\chi(\fu_\lambda)= &\begin{cases}
			v^{\frac{1}{2}\sum_{i\in\I_0}d^\lambda_i +2\{\lambda_f,\lambda_c\} }\frac{a_{2\lambda_f}(v)}{a_{\lambda_f}(v^2)}b_{\lambda_c}(v)
			&\text{if }\lambda\in\cf,
			\\
			0&\text{otherwise,}
		\end{cases}
	\end{align}
	and 
	\begin{align}
		\label{eq:chi3}
		\chi(\fu_\lambda\cdot \K_{\beta^\diamond})= \chi(\fu_\lambda)\cdot \chi(\K_{\beta^\diamond}).
	\end{align}

	We denote by $[0]_v {^{!!}}=[1]_v^{!!}=1$, and for any $k\in \Z_{\ge 1}$,
	\[
	[2k]_v{^{!!}}=[2k]_v[2k-2]_v \cdots [4]_v[2]_v,\qquad [2k-1]_v{^{!!}}=[2k-1]_v[2k-3]_v \cdots [3]_v[1]_v.
	\]
	For convenience, we denote $[-1]_v^{!!}=1$. 
	
	\begin{proposition}
		The linear map $\chi:\widetilde{\mathbf{B}}( Q)\longrightarrow \Q(v^{\frac{1}{2}})$ given in \eqref{eq:chi1}--\eqref{eq:chi3} is $\varrho$-twisted compatible.
	\end{proposition}
	
	\begin{proof}
		Recall that $\widetilde{\mathbf{B}}$ is the Borel subalgebra of $\tUi$ generated by $E_i$ and $K_i$, $i\in\I$. Then the restriction of $\widetilde{\psi}$ to $\widetilde{\mathbf{B}}$ gives the isomorphism $\widetilde{\mathbf{B}}\rightarrow \widetilde{\mathbf{B}}(Q)$ by sending $E_i\mapsto v^{-\frac{1}{2}}\fu_{\alpha_i}$, $K_i\mapsto \K_{\alpha_i^\diamond}$. 
		Using this isomorphism and \cite[Lemma 4.2]{CLPRW}, 
		there exists a unique $\varrho$-twisted compatible map $\chi': \widetilde{\mathbf{B}}(Q) \rightarrow \Q(v^{\frac{1}{2}})$ such that $\chi'(1)=1, \; \chi'(\fu_{\alpha_i})=0, \; \chi'(\K_{i^\diamond})=\varphi(\K_{i^\diamond}, \K_{\varrho i^\diamond})=v^{(i, \varrho i)_Q}$ for $i \in \I$. In particular, we have
		$\chi'(a \cdot b)=\sum\chi'(a_{(1)})\chi'(b_{(2)})\varphi(\varrho (a_{(2)}),b_{(1)})$, for any $a,b\in\widetilde{\mathbf{B}}(Q)$.

		It remains to prove $\chi=\chi'$.
		By definition, 
		we have     \begin{align*}
			\chi'(\fu_\lambda\cdot \K_{\beta^\diamond})= \chi'(\fu_\lambda)\cdot \chi'(\K_{\beta^\diamond}).
		\end{align*}
		So it is enough to prove 
		$\chi(\K_{\beta^\diamond})=\chi'(\K_{\beta^\diamond})$, and $\chi(\fu_\lambda)=\chi'(\fu_\lambda)$, for any $\beta\in\Z^\I$, $\lambda\in\fp$.
		
		Recall that $\beta_{1}, \btau(\beta_1), \beta_{2}, \btau(\beta_2), \dots,\beta_{N_\imath },\btau(\beta_{N_\imath})$ is a $Q$-admissible sequence of $\Phi^+$. For any $\lambda \in \fp$ such that $\lambda(\beta)=0$ if $\beta \in \Phi^+_1$, a similar argument as the proof of Lemma \ref{lem:sqrt-aut} shows that
		\[ a_{\lambda}(v)= \prod_{\beta_k \prec \beta_l \in \Phi^+_0} v^{2 \lambda(\beta_k) \lambda(\beta_l) \dim_{\bfk} \Hom(M_q(\beta_k), M_q(\beta_l)) } \prod_{\beta \in \Phi^+_0} a_{\lambda(\beta)\beta}(v). \]
		By (\ref{function a}) we have
		\begin{align}
			\frac{a_{2\lambda}(v)}{a_{\lambda}(v^2)}= \prod_{\beta_k \prec \beta_l \in \Phi^+_0} v^{4 \lambda(\beta_k) \lambda(\beta_l) \dim_{\bfk} \Hom(M_q(\beta_k), M_q(\beta_l)) } \prod_{\beta \in \Phi^+_0} v^{3\lambda(\beta)^2} (v-v^{-1})^{\lambda(\beta)} [2\lambda(\beta)-1]!!.
			\label{function a(v)/a(q)}
		\end{align}
		
		Now we can prove $\chi=\chi'$. Since $\chi'$ is $\varrho$-twisted compatible, we get
		\[ \chi'(\K_{\beta^\diamond+\alpha_i^\diamond})= \chi'(\K_{\beta^\diamond}) \chi'(\K_{\alpha_i^\diamond}) \varphi(\K_{\varrho \beta^\diamond}, \K_{\alpha_i^\diamond})=v^{(\alpha_i,\varrho(\beta+\alpha_i))_Q} \chi'(\K_{\beta^\diamond}). \]
		An easy induction on $h(\beta):=\sum_{i \in \I} d_i^\beta $ implies $\chi=\chi'$ on $\K_{\beta^\diamond}$ for any $\beta \in \N^\I$. Since
		\[ 1=\chi'(1)=\chi'(\K_{\beta^\diamond} \cdot \K_{-\beta^\diamond})=\chi(\K_{\beta^\diamond})\chi'(\K_{-\beta^\diamond}) v^{-(\varrho \beta,\beta)_Q}, \]
		$\chi=\chi'$ also holds on $\K_{-\beta^\diamond}$ for any $\beta \in \N^\I$. For any $\beta \in \Z^\I$, there exist $\beta_+, \beta_- \in \N^\I$ such that $\beta=\beta_+-\beta_-$. Then
		\[ \chi'(\K_{\beta^\diamond})=\chi(\K_{\beta_+^\diamond})\chi(\K_{-\beta_-^\diamond}) v^{-(\varrho \beta_+,\beta_-)_Q}=v^{\EuFo{\beta, \varrho \beta}+\sum_{i \in \I_0} d^\beta_i}. \]

		For $\lambda \in \fp$, denote $r=\max\{ k \mid \lambda(\beta_k) \neq 0 \text{ or } \lambda(\varrho \beta_k) \neq 0 \}$. Let $\lambda'$ be the function:
		\begin{align*}
			\lambda'(\beta)=\begin{cases}
				\lambda(\beta)&\text{ if }\beta \neq \beta_r,\varrho \beta_r,
				\\
				0&\text{ otherwise}.
			\end{cases}
		\end{align*}
		Let $\mu\in\fp$ be such that 
		$M_q(\lambda)=M_q(\lambda')\oplus M_q(\mu)$. It is clear that each indecomposable direct summand of $M_q(\mu)$ is either $M_q(\beta_r)$ or $M_q(\varrho \beta_r)$. By our assumption, we know $\Ext^1_{\bfk Q}(M_q(\lambda'),M_q(\mu))=0$. In this way
		\[ \fu_\lambda= v^{\dim_\bfk \Hom(M_q(\lambda'), M_q(\mu))} \fu_{\lambda'} \cdot \fu_{\mu}. \] 
		By Lemma \ref{lem: properties of chi}~(1), we have 
		\[\chi'(\fu_{\lambda})=v^{\dim_\bfk \Hom(M_q(\lambda'), M_q(\mu))} \chi'(\fu_{\lambda'})\cdot \chi'(\fu_{\mu}).\]
		The proof is divided into the following two cases.
		
		(i) \underline{$\beta_r \in \Phi^+_0$.} If $\lambda(\beta_r)=1$, then $\chi'(\fu_\mu)=0$. Hence $\chi'(\fu_\lambda)=0$. Of course $\lambda \notin \cf$ in this case. Now suppose $\lambda(\beta_r)\ge 2$. According to Lemma \ref{lem: properties of chi}
		\begin{align*}
			\chi'(\fu_\lambda)=&v^{\dim_\bfk \Hom(M_q(\lambda'), M_q(\mu))+3\lambda(\beta_r)-4}(v-v^{-1})[\lambda(\beta_r)-1]
			\chi'(\fu_{\lambda'}) \chi'(\fu_{(\lambda(\beta_r)-2) \beta_r} \cdot \K_{\beta_r^\diamond}) \\
			=& v^{\dim_\bfk \Hom(M_q(\lambda'), M_q(\mu))+3\lambda(\beta_r)-4}(v-v^{-1})[\lambda(\beta_r)-1]
			\chi'(\fu_{\lambda'})\cdot\chi'(\fu_{(\lambda(\beta_r)-2) \beta_r}) \cdot\chi'(\K_{\beta_r^\diamond}) \\
			=&v^{2\dim_\bfk \Hom(M_q(\lambda'), M_q(\beta_r))+3\lambda(\beta_r)-3+\sum_{i \in \I_0} d^{\beta_r}_i}(v-v^{-1})[\lambda(\beta_r)-1]\cdot
			\chi'(\fu_{\widetilde{\lambda}}),
		\end{align*}
		where $\widetilde{\lambda}\in\fp$ is defined by setting 
		\begin{align*}
			\widetilde{\lambda}(\beta)=\begin{cases}
				\lambda(\beta)&\text{ if }\beta\neq \beta_r,
				\\
				\lambda(\beta_r)-2&\text{ otherwise}.
			\end{cases}
		\end{align*}
		
		On the other hand, if $\lambda \in \cf$, then $2 \mid \lambda(\beta_r)$ and $\widetilde{\lambda} \in \cf$. Using (\ref{function a(v)/a(q)}) we get
		\begin{align*}
			\frac{\chi(\fu_\lambda)}{\chi(\fu_{\widetilde{\lambda}})}
			=& v^{\sum_{i \in \I_0} d^{\beta_r}_i+2 \dim_\bfk \Hom(M_q(\lambda_c),M_q(\beta_r))} \frac{a_{2\lambda_f}(v)}{a_{\lambda_f}(v^2)} \frac{a_{\widetilde{\lambda}}(v^2)}{a_{2\widetilde{\lambda}_f}(v)} \\
			=& v^{\sum_{i \in \I_0} d^{\beta_r}_i+2 \dim_\bfk \Hom(M_q(\lambda'),M_q(\beta_r))+3\lambda(\beta_r)-3}(v-v^{-1}) [\lambda(\beta_r)-1].
		\end{align*}
		
		(ii) \underline{$\beta_r \in \Phi^+_1$.} Without loss of generality we suppose $\lambda(\varrho \beta_r)>0$. By Lemma \ref{lem: properties of chi}~(4), $\chi'(\fu_\lambda)=\chi'(\fu_\mu)=0$ if $\lambda(\beta_r)=0$ (of course $\lambda \notin \cf$ in this case), and if $\lambda(\beta_r)>0$ we have, by Lemma \ref{lem: properties of chi}~(1),
		\begin{align*}
			\chi'(\fu_\lambda)=&v^{\dim_\bfk \Hom(M_q(\lambda'), M_q(\mu))+2\lambda(\beta_r)+\lambda(\varrho \beta_r)-2}(v-v^{-1})[\lambda(\beta_r)]\\ 
			&\cdot\chi'(\fu_{\lambda'}) \cdot\chi'(\fu_{(\lambda(\beta_r)-1)\beta_r+(\lambda(\varrho \beta_r)-1)\varrho \beta_r} \cdot \K_{\beta_r^\diamond}) \\
			=&v^{\dim_\bfk \Hom(M_q(\lambda'), M_q(\mu))+2\lambda(\beta_r)+\lambda(\varrho \beta_r)-2}(v-v^{-1})[\lambda(\beta_r)] \\
			& \cdot\chi'(\fu_{\lambda'})\cdot  \chi'(\fu_{(\lambda(\beta_r)-1)\beta_r+(\lambda(\varrho \beta_r)-1)\varrho \beta_r})\cdot \chi'(\K_{\beta_r^\diamond}) \\
			=&v^{\dim_\bfk \Hom(M_q(\lambda'), M_q(\beta_r) \oplus M_q(\varrho \beta_r))+2\lambda(\beta_r)+\lambda(\varrho\beta_r)-2+\sum_{i \in \I_0} d^{\beta_r}_i}
			(v-v^{-1})[\lambda(\beta_r)] \cdot\chi'(\fu_{\widetilde{\lambda}}),
		\end{align*}
		where $\widetilde{\lambda}\in\fp$ is defined by setting 
		\begin{align*}
			\widetilde{\lambda}(\beta)=\begin{cases}
				\lambda(\beta_r)-1&\text{ if }\beta=\beta_r,
				\\
				\lambda(\varrho\beta_r)-1&\text{ if }\beta=\varrho\beta_r,
				\\
				\lambda(\beta)&\text{ otherwise}.          
			\end{cases}
		\end{align*}
		
		On the other hand, if $\lambda \in \cf$, then $\lambda(\beta_r)=\lambda(\varrho \beta_r)$ and $\widetilde{\lambda} \in \cf$. By Lemma \ref{lem:sqrt-aut} we have
		\begin{align*}
			\frac{\chi(\fu_\lambda)}{\chi(\fu_{\widetilde{\lambda}})}
			=& v^{\sum_{i \in \I_0}d^{\beta_r}_i+2\dim_\bfk \Hom(M_q(\lambda_f),M_q(\beta_r) \oplus M_q(\varrho \beta_r))} \frac{b_{\lambda_c}(v)}{b_{\widetilde{\lambda}_c}(v)} \\
			=& v^{\sum_{i \in \I_0}d^{\beta_r}_i+2\dim_\bfk \Hom(M_q(\lambda_f),M_q(\beta_r) \oplus M_q(\varrho \beta_r))+\dim_\bfk \Hom(M_q(\lambda'_c),M_q(\beta_r) \oplus M_q(\varrho \beta_r))} \\
			& \cdot v^{3\lambda(\beta_r)-2} (v-v^{-1}) [\lambda(\beta_r)] \\
			=& v^{\sum_{i \in \I_0}d^{\beta_r}_i+\dim_\bfk \Hom(M_q(\lambda'),M_q(\beta_r) \oplus M_q(\varrho \beta_r))+3\lambda(\beta_r)-2}(v-v^{-1}) [\lambda(\beta_r)] .
		\end{align*}
		
		We use an induction on $\wt(\lambda):=\sum_{\beta \in \Phi^+} \lambda(\beta)$. In both cases $\wt(\widetilde{\lambda}) < \wt(\lambda)$ and $\lambda \in \cf$ if and only if $\widetilde{\lambda} \in \cf$. Therefore, $\chi'(\fu_\lambda)=0=\chi(\fu_\lambda)$ if $\lambda \notin \cf$. Now suppose $\lambda \in \cf$, then the above argument shows that
		\[ \chi'(\fu_\lambda)=\frac{\chi(\fu_\lambda)}{\chi(\fu_{\widetilde{\lambda}})} \chi'(\fu_{\widetilde{\lambda}}). \]
		The induction on $\wt(\lambda)$ completes the proof.
	\end{proof}

	With the help of $\varrho$-twisted compatible linear map $\chi$ defined in \eqref{eq:chi1}--\eqref{eq:chi3}, we can obtain the following theorem, which explains the natural embedding $\imath:\tUi\to \tU$ in Hall algebra setting.
	
	\begin{theorem}
		\label{thm:embedding}
		There exists an injective $\Q(v^{\frac12})$-algebra homomorphism
		\begin{align}
			\label{eq:embedding}
			\widetilde{\Omega}: \widetilde{\ch}(Q,\varrho)&\longrightarrow \widetilde{\ch}(Q^{\dbl},\swa)
			\\\notag
			\fu_{\lambda} &\mapsto \sum_{\lambda_1,\lambda_2 \in \fp,\nu\in\cf} v^{\Xi} f^\lambda_{\lambda_1,\nu,\lambda_2}(v) 
			\cdot\frac{a_{2\nu_f}(v)}{a_{\nu_f}(v^2)}b_{\nu_c}(v) \cdot
			\fu_{\varrho \lambda_2} \oplus \fu_{ \lambda_1^\diamond} *\K_{\nu+\lambda_2},
			\\\notag
			\K_{\alpha} &\mapsto v^{\sum_{i \in \I_0} d_i^\alpha} \K_{\varrho \alpha} * \K_{\alpha^\diamond}, 
		\end{align}
		where
		\begin{align*} \Xi=& \EuFo{\lambda_1-\varrho \lambda_2,\nu+\lambda_2}_Q
			+\sum_{i\in\I_0}d_i^{\lambda_2}+\frac{1}{2}\sum_{i\in\I_0}d_i^{\nu}+2\{\nu_f,\nu_c\}.\end{align*}
		Moreover, we have the following commutative diagram:
		$$ \xymatrix{ \tUi \ar[rr]^{\imath }\ar[d]^{\widetilde{\psi}} && \tU \ar[d]^{\widetilde{\psi}^{\dbl}} 
			\\
			\widetilde{\ch}(Q,\varrho) \ar[rr]^{\widetilde{\Omega}} && \widetilde{\ch}(Q^{\dbl},\swa)} $$
	\end{theorem}

	\begin{proof}
		By \cite[Proposition 2.12]{CLPRW} there exists an algebra homomorphism $\widetilde{\xi}_\varrho: \widetilde{\mathbf{B}}^\imath_\varrho \rightarrow D(\widetilde{\mathbf{B}}(Q))$ defined by
		\[ a \mapsto \sum \chi(a_{(2)}) \cdot \varrho(a_{(3)}) \otimes a_{(1)}, \quad \forall a \in \widetilde{\mathbf{B}}^\imath_\varrho.\]
		We claim $\widetilde{\Omega}=(\widetilde{\Phi}^{\dbl})^{-1} \circ \widetilde{\xi}_\varrho \circ \widetilde{\Phi}$. Indeed, for any $\lambda \in \fp$ and $\alpha \in \Z^\I$ we have
		\begin{align*}
			&(\widetilde{\Phi}^{\dbl})^{-1} \circ \widetilde{\xi}_\varrho \circ \widetilde{\Phi}(\fu_{ \lambda} * \K_{\alpha})= v^{(\lambda,\alpha)_Q+\EuFo{\alpha,\varrho \alpha}_Q} (\widetilde{\Phi}^{\dbl})^{-1} \widetilde{\xi}_\varrho (\fu_{ \lambda} \cdot \K_{\varrho \alpha^\diamond}) \\
			=& v^{(\lambda,\alpha)_Q+\EuFo{\alpha,\varrho \alpha}_Q}
			(\widetilde{\Phi}^{\dbl})^{-1}\Bigg(\sum_{\lambda_1,\nu,\lambda_2 \in \fp}v^{\EuFo{\lambda_1, \nu+ \lambda_2}_Q+\EuFo{ \nu,\lambda_2}_Q} f^{ \lambda}_{\lambda_1,  \nu, \lambda_2}(v) \cdot \chi(\fu_{ \nu}\cdot \K_{\varrho\alpha^\diamond+ \lambda_2^\diamond}) \\
			& \qquad \qquad\qquad\qquad \qquad\cdot \varrho(\fu_{ \lambda_2} \cdot \K_{\varrho \alpha^\diamond}) \otimes \fu_{ \lambda_1} \cdot \K_{\varrho\alpha^\diamond+ \nu^\diamond+\lambda_2^\diamond} \Bigg) \\
			=& \sum_{\lambda_1,\lambda_2 \in \fp,\nu \in \cf} v^{\EuFo{\lambda_1-\varrho \lambda_2,\nu+\lambda_2}_Q+\sum_{i \in I_0} d^{\lambda_2}_i+\frac{1}{2}\sum_{i \in \I_0} d^\nu_i+2\{\nu_f,\nu_c\}+\sum_{i \in \I_0}d_i^{ \alpha}}
			f^\lambda_{\lambda_1,\nu,\lambda_2}(v)\cdot \frac{a_{2\nu_f}(v)}{a_{\lambda_f}(v^2)}b_{\nu_c}(v) \\
			&\qquad \qquad\cdot 
			\fu_{\varrho \lambda_2} \oplus \fu_{ \lambda_1^\diamond} *\K_{\varrho\alpha+\nu+\lambda_2}*\K_{ \alpha^\diamond}.
		\end{align*}
		In particular, for any simple root $\alpha_i$ we get
		\[ \widetilde{\Omega}(\fu_{\alpha_i})= \fu_{\alpha_i^\diamond}+ \fu_{\varrho \alpha_i} * \K_{ \alpha_i}, \quad  \widetilde{\Omega}(\K_{\alpha_i})=v^{\EuFo{\alpha_i,\varrho \alpha_i}_Q} \K_{\varrho \alpha_i}*\K_{\alpha_i^\diamond}. \]
		Thus, for any $i \in \I$
		\begin{align*} \widetilde{\psi}^{\dbl}\circ \imath(B_i)&= \widetilde{\psi}^{\dbl}(F_i+E_{\varrho i} K_i')=v^{-\frac{1}{2}}(\fu_{\alpha_i^\diamond}+\fu_{\varrho \alpha_i}*\K_{\alpha_i}), \\
			\widetilde{\psi}^{\dbl} \circ \imath(\widetilde{k_i})&=\widetilde{\psi}^{\dbl}(K_iK_{\varrho i}')=\K_{\alpha_i^\diamond} * \K_{\varrho \alpha_i}, \\
			\widetilde{\Omega} \circ \widetilde{ \psi}(B_i)&=v^{-\frac{1}{2}} \widetilde{\Omega}(\fu_{\alpha_i})=v^{-\frac{1}{2}}(\fu_{\alpha_i^\diamond}+\fu_{\varrho \alpha_i}*\K_{\alpha_i}), \\
			\widetilde{\Omega} \circ \widetilde{ \psi}(\widetilde{k_i})&=v^{-\EuFo{\alpha_i,\varrho \alpha_i}}\widetilde{\Omega} (\K_{\alpha_i})= \K_{\varrho \alpha_i}*\K_{\alpha_i^\diamond}.\end{align*}
		Therefore, $\widetilde{\psi}^{\dbl}\circ \imath=\widetilde{\Omega} \circ \widetilde{\psi}$ holds on the generators of $\tUi$. Since $\widetilde{\psi},\widetilde{\psi}^{\dbl}$ are isomorphisms, we can see that $\widetilde{\Omega}$ is injective, which completes the proof.
	\end{proof}
	
	For the case $\varrho=\Id$, the split case, we have $\Phi^+_0=\Phi^+$. Then $\cf=\{\lambda\in\fp\mid \lambda(\beta)\in2\N,\forall \beta\in\Phi^+\}$. And for any $\lambda\in\cf$, we see $\lambda_f=\frac{\lambda}{2}$. 
	As a corollary of Theorem \ref{thm:embedding}, we have the following.
	\begin{corollary}
		If $\varrho=\Id$, then there exists an injective $\Q(v^{\frac12})$-algebra homomorphism:
		\begin{align}
			\widetilde{\Omega}: \widetilde{\ch}(Q,\varrho)&\longrightarrow \widetilde{\ch}(Q^{\dbl},\swa)
			\\\notag
			\fu_{\lambda} &\mapsto \sum_{\lambda_1,\lambda_2 \in \fp,\nu\in\cf} v^{\Xi} f^\lambda_{\lambda_1,\nu,\lambda_2}(v) 
			\cdot\frac{a_{2\nu_f}(v)}{a_{\nu_f}(v^2)} \cdot
			\fu_{\lambda_2} \oplus \fu_{ \lambda_1^\diamond} *\K_{\nu+\lambda_2},
			\\\notag
			\K_{\alpha} &\mapsto v^{\sum_{i \in \I} d_i^\alpha} \K_{ \alpha} * \K_{\alpha^\diamond}, 
		\end{align}
		where
		\begin{align*} \Xi=& \EuFo{\lambda_1- \lambda_2,\nu+\lambda_2}_Q
			+\sum_{i\in\I}d_i^{\lambda_2}+\frac{1}{2}\sum_{i\in\I}d_i^{\nu}.\end{align*}    
	\end{corollary}
	
	\begin{remark}
		We expect that the result in Theorem \ref{thm:embedding} also holds for arbitrary quivers $Q$ (not necessarily of Dynkin type). However, the $\varrho$-twisted compatible linear map $\chi:\widetilde{\mathbf{B}}( Q)\rightarrow \Q(v^{\frac12})$ in this general setting is much more complicated than the one formulated in \eqref{eq:chi1}--\eqref{eq:chi3}.
	\end{remark}  
	
	At the end of this section, we consider the tensor algebra $\widetilde{\ch}(Q,\varrho)\otimes\widetilde{\ch}(Q^{\dbl},\swa)$. We have the following result, which explains the coideal algebra structure of $\tUi$ in $\tU$.

	\begin{theorem}
		\label{thm:coproduct}
		There exists an $\Q(v^{\frac12})$-algebra homomorphism
		\begin{align}
			\label{eq:coproduct}
			&\widetilde{\Delta}: \widetilde{\ch}(Q,\varrho)\longrightarrow \widetilde{\ch}(Q,\varrho)\otimes\widetilde{\ch}(Q^{\dbl},\swa)
			\\\notag
			\fu_\lambda&\mapsto \sum_{\lambda_1,\lambda_2,\lambda_3,\nu \in \fp}
			v^{\Xi'} f^\lambda_{\lambda_1,\nu,\lambda_2,\varrho \nu,\lambda_3}(v) \cdot a_{\nu}(v) \cdot
			(\fu_{\lambda_2} * \K_{\varrho \lambda_3+\nu} )\otimes (\fu_{\varrho \lambda_3}\oplus \fu_{\lambda_1^\diamond} * \K_{\lambda_2+\lambda_3+\nu+\varrho \nu}),
			\\\notag
			\K_\alpha&\mapsto\K_{\alpha} \otimes (\K_{\varrho \alpha}*\K_{\alpha^\diamond}),
		\end{align}
		where
		\[ \Xi'= \EuFo{\lambda_1,\lambda_2+\lambda_3}_Q+\EuFo{\lambda_2,\lambda_3}_Q+\EuFo{\lambda_3,\varrho \lambda_3}_Q-(\lambda_2+\lambda_3,\varrho \lambda_3)_Q+\EuFo{\lambda_1-\lambda_3,\nu+\varrho \nu}_Q-\EuFo{\lambda_2,\nu-\varrho \nu}_Q.\]
		Moreover, we have the following commutative diagram:
		\[ \xymatrix{ \tUi \ar[rr]^{\Delta\circ\imath }\ar[d]^{\widetilde{\psi}} && \tUi\otimes\tU \ar[d]^{\widetilde{\psi}\otimes\widetilde{\psi}^{\dbl}} 
			\\
			\widetilde{\ch}(Q,\varrho) \ar[rr]^{\widetilde{\Delta}\qquad\qquad} && \widetilde{\ch}(Q,\varrho)\otimes\widetilde{\ch}(Q^{\dbl},\swa)} \]
	\end{theorem}
	
	\begin{proof}
		By \cite[Proposition 2.13]{CLPRW} there exists an algebra homomorphism $\widetilde{\Psi}_\varrho: \widetilde{\mathbf{B}}_\varrho^\imath \rightarrow \widetilde{\mathbf{B}}_\varrho^\imath \otimes D(\widetilde{\mathbf{B}}(Q))$ given by
		\[ a \mapsto \sum \varphi(\varrho(a_{(4)}), a_{(2)}) \cdot a_{(3)} \otimes (\varrho(a_{(5)}) \otimes a_{(1)}). \]
		We claim that $\widetilde{\Delta}=(\widetilde{\Phi} \otimes \widetilde{\Phi}^{\dbl})^{-1} \circ \widetilde{\Psi}_\varrho \circ \widetilde{\Phi}$. Indeed, for any $\lambda \in \fp$, a calculation shows that
		\begin{align*}
			&\widetilde{\Psi}_\varrho \circ \widetilde{\Phi}(\fu_\lambda)
			=\widetilde{\Psi}_\varrho(\fu_{ \lambda}) \\
			=& \sum_{\lambda_1,\lambda_2,\lambda_3,\nu \in \fp}
			v^{\EuFo{\lambda_1,\lambda_2+\lambda_3+\nu+\varrho \nu}_Q+\EuFo{\nu,\lambda_2+\lambda_3+\varrho \nu}_Q+\EuFo{\lambda_2,\lambda_3+\varrho \nu}_Q+\EuFo{\varrho \nu, \lambda_3}_Q+(\varrho \lambda_3, \lambda_2+\lambda_3+\varrho \nu)_Q} \\
			& \qquad \qquad \cdot f^\lambda_{\lambda_1,\nu,\lambda_2,\varrho \nu,\lambda_3}(v) \cdot a_{\nu}(v) \cdot
			(\fu_{\lambda_2} \cdot \K_{\lambda_3^\diamond+\varrho\nu^\diamond} )\otimes (\fu_{\varrho \lambda_3}\otimes \fu_{\lambda_1} \cdot \K_{\lambda_2^\diamond+\lambda_3^\diamond+\nu^\diamond+\varrho \nu^\diamond}).
		\end{align*}
		Applying $(\widetilde{\Phi} \otimes \widetilde{\Phi}^{\dbl})^{-1}$ on the above formula one gets
		\begin{align*}
			&(\widetilde{\Phi} \otimes \widetilde{\Phi}^{\dbl})^{-1} \circ \widetilde{\Psi}_\varrho \circ \widetilde{\Phi}(\fu_\lambda) \\
			=&\sum_{\lambda_1,\lambda_2,\lambda_3,\nu \in \fp}
			v^{\EuFo{\lambda_1,\lambda_2+\lambda_3+\nu+\varrho \nu}_Q+\EuFo{\nu,\lambda_2+\lambda_3+\varrho \nu}_Q+\EuFo{\lambda_2,\lambda_3+\varrho \nu}_Q+\EuFo{\varrho \nu, \lambda_3}_Q+(\varrho \lambda_3, \lambda_2+\lambda_3+\varrho \nu)_Q} \\
			&\qquad \qquad \cdot v^{-(\lambda_2,\varrho \lambda_3+\nu)_Q-\EuFo{\lambda_3+\varrho \nu, \varrho \lambda_3+\nu}_Q-(\varrho \lambda_3,\lambda_2+\lambda_3+\nu+\varrho \nu)_Q}\\
			&\qquad \qquad \cdot f^\lambda_{\lambda_1,\nu,\lambda_2,\varrho \nu,\lambda_3}(v) \cdot a_{\nu}(v) \cdot
			(\fu_{\lambda_2} * \K_{\varrho \lambda_3+\nu} )\otimes (\fu_{\varrho \lambda_3}\oplus \fu_{\lambda_1^\diamond} * \K_{\lambda_2+\lambda_3+\nu+\varrho \nu}) \\
			=&\sum_{\lambda_1,\lambda_2,\lambda_3,\nu \in \fp}
			v^{\EuFo{\lambda_1,\lambda_2+\lambda_3}_Q+\EuFo{\lambda_2,\lambda_3}_Q+\EuFo{\lambda_3,\varrho \lambda_3}_Q-(\lambda_2+\lambda_3,\varrho \lambda_3)_Q+\EuFo{\lambda_1-\lambda_3,\nu+\varrho \nu}_Q-\EuFo{\lambda_2,\nu-\varrho \nu}_Q}\\
			& \qquad \qquad \cdot f^\lambda_{\lambda_1,\nu,\lambda_2,\varrho \nu,\lambda_3}(v) \cdot a_{\nu}(v) \cdot
			(\fu_{\lambda_2} * \K_{\varrho \lambda_3+\nu} )\otimes (\fu_{\varrho \lambda_3}\oplus \fu_{\lambda_1^\diamond} * \K_{\lambda_2+\lambda_3+\nu+\varrho \nu}).
		\end{align*}
		For any $\alpha \in \Z^\I$ we have
		\begin{align*}
			(\widetilde{\Phi} \otimes \widetilde{\Phi}^{\dbl})^{-1} \circ \widetilde{\Psi}_\varrho \circ \widetilde{\Phi}(\K_{\alpha})=&v^{\EuFo{\alpha,\varrho \alpha}_Q}(\widetilde{\Phi} \otimes \widetilde{\Phi}^{\dbl})^{-1} \circ \widetilde{\Psi}_\varrho (\K_{\varrho \alpha^\diamond}) \\
			=&v^{\EuFo{\alpha,\varrho \alpha}_Q+(\alpha,\varrho \alpha)_Q} (\widetilde{\Phi} \otimes \widetilde{\Phi}^{\dbl})^{-1} (\K_{\varrho \alpha^\diamond} \otimes (\K_{\alpha^\diamond} \otimes \K_{\varrho \alpha^\diamond})) \\
			=& \K_{\alpha} \otimes (\K_{\varrho \alpha}*\K_{\alpha^\diamond}).
		\end{align*}
		In particular, for any simple root $\alpha_i$,
		\[ \widetilde{\Delta}(\fu_{\alpha_i})=1 \otimes \fu_{\alpha_i^\diamond}+\fu_{\alpha_i} \otimes \K_{\alpha_i}+v^{-\EuFo{\alpha_i,\varrho \alpha_i}_Q} \K_{\varrho \alpha_i} \otimes (\fu_{\varrho \alpha_i}*\K_{\alpha_i}).\]
		Thus, for any $i \in \I$,
		\begin{align*}
			(\widetilde{\psi} \otimes \widetilde{ \psi}^{\dbl}) \circ \Delta \circ \imath(B_i)=&(\widetilde{\psi} \otimes \widetilde{ \psi}^{\dbl})(1 \otimes F_i+B_i \otimes K_i'+\widetilde{k}_{\varrho i} \otimes E_{\varrho i}K_i')\\
			=&v^{-\frac{1}{2}}(1 \otimes \fu_{\alpha_i^\diamond}+\fu_{\alpha_i} \otimes \K_{\alpha_i}+v^{-\EuFo{\alpha_i,\varrho \alpha_i}_Q} \K_{\varrho \alpha_i} \otimes \fu_{\varrho \alpha_i}*\K_{\alpha_i}),
			\\ \widetilde{\Delta} \circ \widetilde{ \psi}(B_i)=& v^{-\frac{1}{2}}\widetilde{\Delta}(\fu_{\alpha_i})=v^{-\frac{1}{2}}(1 \otimes \fu_{\alpha_i^\diamond}+\fu_{\alpha_i} \otimes \K_{\alpha_i}+v^{-\EuFo{\alpha_i,\varrho \alpha_i}_Q} \K_{\varrho \alpha_i} \otimes \fu_{\varrho \alpha_i}*\K_{\alpha_i}),\\
			(\widetilde{\psi} \otimes \widetilde{ \psi}^{\dbl}) \circ \Delta \circ \imath(\widetilde{k}_i)=&(\widetilde{\psi} \otimes \widetilde{ \psi}^{\dbl})(\widetilde{k}_i \otimes K_iK_{\varrho i}')=v^{-\EuFo{\alpha_i,\varrho \alpha_i}_Q} \K_{\alpha_i} \otimes (\K_{\varrho \alpha_i}*\K_{\alpha_i^\diamond}), \\
			\widetilde{\Delta} \circ \widetilde{ \psi}(\widetilde{k}_i)=&v^{-\EuFo{\alpha_i,\varrho \alpha_i}_Q} \widetilde{\Delta}(\K_{\alpha_i})=v^{-\EuFo{\alpha_i,\varrho \alpha_i}_Q} \K_{\alpha_i} \otimes (\K_{\varrho \alpha_i}*\K_{\alpha_i^\diamond}).\end{align*}
		Therefore, $(\widetilde{\psi} \otimes \widetilde{ \psi}^{\dbl}) \circ \Delta \circ \imath=\widetilde{\Delta} \circ \widetilde{ \psi}$ holds on the generators of $\tUi$.
	\end{proof}


	\section{Integral forms and dual canonical bases}
	\label{sec:dCB}
	
	In \cite{LP25}, Pan and the first author introduce integral forms on $\tU$, $\tUi$, and construct their dual canonical bases. 
	In this section, we shall prove that the integral forms of $\tU$ and $\tUi$ are compatible under the embedding and coproduct.

	\subsection{Integral forms}
	
	Let $\cz=\Z[v^{\frac12},v^{-\frac12}]$ and $\cz_\sqq=\Z[\sqq^{\frac12},\sqq^{-\frac12}]$. We define the integral form $\widetilde{\ch}(\bfk Q,\varrho)_{\cz_\sqq}$ of $\widetilde{\ch}(\bfk Q,\varrho)$ to be the $\cz$-module spanned by $\{[\K_\alpha]\ast[M_q(\mu)]\mid \alpha\in\Z^\I,\mu\in\fp\}$. Note that $\widetilde{\ch}(\bfk Q,\varrho)_{\cz_\sqq}$ is a free $\cz_\sqq$-module. 
	Using \eqref{eq:multiplication}, we know that $\widetilde{\ch}(\bfk Q,\varrho)_{\cz_\sqq}$ is a $\cz_\sqq$-algebra, and $\widetilde{\ch}(\bfk Q,\varrho)_{\cz_\sqq}\otimes \Q(\sqq^{\frac12})\cong \widetilde{\ch}(\bfk Q,\varrho)$. Similarly, we can define the (generic) integral form $\widetilde{\ch}(Q,\varrho)_\cz$.

	Using the  isomorphism of $\Q(v^{1/2})$-algebras
	$\widetilde{\psi}:\tUi\rightarrow \widetilde{\ch}(Q,\btau)$, 
	we can define the integral form of $\tUi$ as the preimage of $\widetilde{\ch}(Q,\varrho)_\cz$, which is denoted by $\tUi_\cz$. It is also a  free $\cz$-module. 
	By \cite[Corollary 7.10]{LP25}, we know $\tUi_\cz$ is independent of the orientation of $Q$.

	Considering the iquiver $( Q^{\dbl},\swa)$, we can get the integral forms $\widetilde{\ch}(\bfk Q^{\dbl},\swa)_{\cz_\sqq}$, $\widetilde{\ch}( Q,\swa)_{\cz}$, and 
	$\tU_\cz$. Let $\widetilde{\mathbf{B}}(\bfk Q)_{\cz_v}$ be the free $\cz_\sqq$-submodule  of $\widetilde{\ch}(\bfk Q^{\dbl},\swa)_{\cz_\sqq}$ spanned by $\{[\K_{\alpha^\diamond}]*[M_q(\mu)]\mid \alpha\in\Z^\I,\mu\in\fp\}$. Then $\widetilde{\mathbf{B}}(\bfk Q)_{\cz_\sqq}$ is a $\cz$-subalgebra of $\widetilde{\ch}(\bfk Q^{\dbl},\swa)_{\cz_\sqq}$. 
	Similarly, we can define $\widetilde{\mathbf{B}}( Q)_\cz$. 
	
	\begin{lemma}
		\label{lem:laurent-chi}
		We have $a_\lambda(v^2)\mid a_{2\lambda}(v)$ for any $\lambda\in\fp$. Moreover, we have $\chi: \widetilde{\mathbf{B}}( Q)_\cz\rightarrow\cz$ induced by the linear map 
		$\chi: \widetilde{\mathbf{B}}( Q)\rightarrow \Q(v^{\frac{1}{2}})$. 
	\end{lemma}
	
	\begin{proof}
		
		Let us prove the first statement. 
		Fix a $Q$-admissible sequence $\beta_1,\dots,\beta_N$ of $\Phi^+$. Then $\Hom(M_q(\beta_i),M_q(\beta_j))\neq0$ only if $i\leq j$. Denote $r_i=\lambda(\beta_i)$ for $1\leq i\leq N$. 
		Note that $M_q(\lambda)=\oplus_{i=1}^NM_q(r_i\beta_i)$. Then
		\begin{align*}
			a_{M_q(\lambda)}=q^{\sum_{i<j}r_ir_j\dim_\bfk\Hom(M_q(\beta_i),M_q(\beta_j)) }\prod_{\beta\in\Phi^+}a_{M_q(r_i\beta_i)}.
		\end{align*}
		The formulas for $a_{M_{q^2}(\lambda)}$ and $a_{M_{q}(2\lambda)}$ are similar. It is enough to prove $a_{m\beta}(v^2)\mid a_{2m\beta}(v)$ for any $\beta\in\Phi^+$, $m\in\N$. Note that \begin{align*}
			|\Aut(M_q(2m\beta))|=&(q^{2m}-1)(q^{2m}-q)\cdots(q^{2m}-q^{2m-1}),\\ 
			|\Aut(M_{q^2}(m\beta))|=&(q^{2m}-1)(q^{2m}-q^2)\cdots(q^{2m}-q^{2(m-1)}).
		\end{align*}
		One can see 
		\begin{align*}
			\frac{|\Aut(M_q(2m\beta))|}{|\Aut(M_{q^2}(m\beta))|}=(q^{2m}-q)(q^{2m}-q^3)\cdots (q^{2m}-q^{2m-1}).
		\end{align*}
		Since the finite field $\F_q$ is arbitrary, we have proved $a_{m\beta}(v^2)\mid a_{2m\beta}(v)$, and then the desired result follows.
		
		The second statement follows by the first one and its definition; see \eqref{eq:chi1}--\eqref{eq:chi3}.
	\end{proof}
	
	With the help of Lemma \ref{lem:laurent-chi} and Theorems \ref{thm:embedding}--\ref{thm:coproduct}, we can obtain the following proposition.
	
	\begin{proposition}
		\label{prop:interal-embedding-coproduct}
		\begin{enumerate}
			\item The $\Q(v^{\frac12})$-algebra homomorphism 
			$\widetilde{\Omega}: \widetilde{\ch}(Q,\varrho)\rightarrow \widetilde{\ch}(Q^{\dbl},\swa)$ 
			induces  a $\cz$-algebra homomorphism
			$\widetilde{\Omega}: \widetilde{\ch}(Q,\varrho)_\cz\rightarrow \widetilde{\ch}(Q^{\dbl},\swa)_\cz$. 
			\item 
			The $\Q(v^{\frac12})$-algebra homomorphism
			$\widetilde{\Delta}: \widetilde{\ch}(Q,\varrho)\rightarrow \widetilde{\ch}(Q,\varrho)\otimes\widetilde{\ch}(Q^{\dbl},\swa)$ 
			induces a $\cz$-algebra homomorphism
			$\widetilde{\Delta}: \widetilde{\ch}(Q,\varrho)_\cz\rightarrow \widetilde{\ch}(Q,\varrho)_\cz\otimes\widetilde{\ch}(Q^{\dbl},\swa)_\cz$. 
		\end{enumerate}
	\end{proposition}

	\begin{proof}
		It is enough to prove that all the coefficients appearing in \eqref{eq:embedding}, \eqref{eq:coproduct} are in $\cz$, which follows from \eqref{eq:Hall-poly}, Lemma \ref{lem:sqrt-aut}, Lemma \ref{lem:laurent-chi}. 
	\end{proof}
	
	As a corollary of Proposition \ref{prop:interal-embedding-coproduct}, by using the definition of $\tUi_\cz$ and $\tU_\cz$, we have the following.
	\begin{corollary}
		\begin{enumerate}
			\item The natural embedding satisfies $\imath(\tUi_\cz)\subseteq\tU_\cz$.
			\item The coproduct satisfies $\Delta(\tUi_\cz)\subseteq \tUi_\cz\otimes\tU_\cz$.
		\end{enumerate}
	\end{corollary}

	\subsection{Dual canonical bases}
	
	With the help of the isomorphism $\widetilde{\psi}: \tUi\rightarrow \widetilde{\ch}(Q,\varrho)$, we can transform the bar-involution of $\tUi$ to $\widetilde{\ch}(Q,\varrho)$. More explicitly,
	\[\ov{\fu_{\alpha_i}}=v^{-1}\fu_{\alpha_i},\quad \ov{\K_{\alpha_i}}=\K_{\alpha_i}.\]
	
	\begin{lemma}[{\cite[Lemma 5.16]{LP25}}]\label{lem:bar on iHall integral}
		The bar-involution of $\widetilde{\ch}(Q,\varrho)$ preserves $\widetilde{\ch}(Q,\btau)_\cz$. Moreover, the bar-involution of $\tUi$ preserves $\tUi_\cz$.
	\end{lemma}
	
	For $\lambda\in\mathfrak{P}$, we define
	\begin{equation}\label{eq:HA element H_lambda}
		\mathfrak{U}_\lambda=v^{-\dim\End_{\bfk Q}(M_q(\lambda))+\frac{1}{2}\langle M_q(\lambda),M_q(\lambda)\rangle_Q}\fu_\lambda.
	\end{equation}
	
	Following \cite{BG17a,LP25}, we define 
	\[\K_\alpha\diamond\fu_\lambda=v^{\frac{1}{2}(\alpha-\varrho\alpha,\,\dimv{M_q(\lambda)})_Q}\K_\alpha\ast\fu_\lambda.\]
	For any $\alpha\in\N^\I$ and $\lambda\in\mathfrak{P}$,
	\begin{equation}\label{eq:diamond action and bar}
		\ov{\K_\alpha\diamond \fu_\lambda}=\K_\alpha\diamond\ov{\fu_\lambda}.
	\end{equation}
	
	Let $\prec$ be the partial order on $\mathfrak{P}$ defined by orbit closure: we say that $\lambda\prec\mu$ if the orbit $\mathfrak{O}_{M_q(\lambda)}$ is contained in the closure of $\mathfrak{O}_{M_q(\mu)}$; see \cite{Lus90} or \cite[\S1.6]{DDPW}. We define a partial order on $\Z^\I\times\mathfrak{P}$: we say $(\alpha,\lambda)\prec(\beta,\mu)$ if $\alpha+\btau(\alpha)+\dimv M_q(\lambda)=\beta+\btau(\beta)+\dimv M_q(\mu)$ and either $\alpha\prec\beta$ (i.e. $\alpha\neq\beta$ and $\beta-\alpha\in\N^\I$) or $\alpha=\beta$ and $\lambda\prec\mu$.

	\begin{theorem}[{\cite[Theorem 5.18]{LP25}}]\label{iHA dCB theorem}
		For each $\alpha\in\Z^\I$ and $\lambda\in\mathfrak{P}$, there exists a unique element $\mathfrak{L}_{\alpha,\lambda}\in\widetilde{\ch}(Q,\varrho)$ such that $\ov{\mathfrak{L}_{\alpha,\lambda}}=\mathfrak{L}_{\alpha,\lambda}$ and
		\[
		\mathfrak{L}_{\alpha,\lambda}-\K_\alpha\diamond \mathfrak{U}_\lambda\in\sum_{(\beta,\mu)}v^{-1}\Z[v^{-1}]\cdot \K_\beta\diamond \mathfrak{U}_\mu.
		\]
		Moreover, $\mathfrak{L}_{\alpha,\lambda}$ satisfies 
		\[
		\mathfrak{L}_{\alpha,\lambda}-\K_\alpha\diamond \mathfrak{U}_\lambda\in\sum_{(\alpha,\lambda)\prec(\beta,\mu)}v^{-1}\Z[v^{-1}]\cdot \K_\beta\diamond \mathfrak{U}_\mu,
		\]
		and $\mathfrak{L}_{\alpha,\lambda}=\K_\alpha\diamond \mathfrak{L}_{0,\lambda}$.
	\end{theorem}
	
	Then $\mathbf{C}^\imath=\{\K_\alpha\diamond \mathfrak{L}_\lambda\mid \alpha\in\Z^\I,\lambda\in\mathfrak{P}\}$ is a $\cz$-basis of $\widetilde{\ch}(Q,\varrho)_\cz$, called the dual canonical basis. It is proved in \cite{LP25} that $\mathbf{C}^\imath$ does not depend on the orientation of $Q$. By transforming to $\tUi$ via $\widetilde{\psi}$, we can get a $\cz$-basis  of $\tUi_\cz$, also denoted by $\mathbf{C}^\imath=\{\K_\alpha\diamond \mathfrak{L}_\lambda\mid \alpha\in\Z^\I,\lambda\in\mathfrak{P}\}$.
	Considering $\widetilde{\ch}(Q^{\dbl},\swa)$, we can get the dual canonical bases of $\widetilde{\ch}(Q^{\dbl},\swa)_\cz$, and $\tU_\cz$. Both of them are denoted by $\mathbf{C}=\{\K_\alpha\diamond \K_{\beta^\diamond}\diamond \mathfrak{L}_{\lambda,\mu^\diamond}\mid \alpha,\beta\in\Z^\I,\lambda,\mu\in\fp\}$. It is remarkable that the structure constants of the multiplication are in $\N[v^{\frac12},v^{-\frac12}]$; see \cite{LP25}.  
	
	Then Proposition \ref{prop:interal-embedding-coproduct}  gives the following result.
	\begin{proposition}
		\label{prop:integral-dCB}
		For any $b\in \mathbf{C}^\imath$, we have
		\begin{align}
			\imath(b)=&\sum_{c\in\mathbf{C}}f_c(v)c,\qquad
			\widetilde{\Delta}(b)=  \sum_{b'\in\mathbf{C}^\imath,c\in\mathbf{C}}g_{b',c''}(v)b'\otimes c,
			\\
			&(\imath\otimes1)\circ \widetilde{\Delta}(b)=   \sum_{c\in\mathbf{C},c'\in\mathbf{C}}h_{c,c'}(v)c\otimes c',
		\end{align}
		where all the constants $f_c(v), g_{b',c}(v),h_{c,c'}(v)\in\cz$. 
	\end{proposition}

	
	\begin{conjecture}[{\cite[Conjectures 1.7.2, 1.7.3]{LP25b}}]
		\label{conj:positive}
		For any $b\in\mathbf{C}^\imath$, all the constants $f_c(v), g_{b',c}(v),h_{c,c'}(v)\in\N[v^{\frac12},v^{-\frac12}]$. In particular, for $\tU$, the structure constants of the coproduct are in $\N[v^{\frac12},v^{-\frac12}]$.
	\end{conjecture}
	
	This conjecture has been proved for $\mathfrak{sl}_2$ in \cite{CZ26} by direct computations. We know $\fU_\beta\in \mathbf{C}^\imath$ for any $\beta\in\Phi^+$; see \cite[Corollary 5.22]{LP25}. At the end of this section, we compute $\widetilde{\Omega}(\fU_\beta)$ for type $A$. 
	We first give the formula of $\widetilde{\Omega}(\fU_{\lambda})$ for general $\widetilde{\ch}(Q,\varrho)$ and $\lambda \in \fp$:
	\begin{align*}
		\widetilde{\Omega}(\fU_{\lambda})=&v^{-\dim_\bfk \End(M_q(\lambda))+\frac{1}{2}\EuFo{\lambda,\lambda}_Q} \widetilde{\Omega}(\fu_{\lambda})\\
		=&v^{-\dim_\bfk \End(M_q(\lambda))+\frac{1}{2}\EuFo{\lambda,\lambda}_Q}\sum_{\lambda_1,\lambda_2 \in \fp, \nu \in \cf} v^{\EuFo{\lambda_1-\varrho \lambda_2,\nu+\lambda_2}_Q+\sum_{i \in \I_0} d^{\lambda_2}_i+\frac{1}{2}\sum_{i \in \I_0} d^\nu_i+2\{\nu_f,\nu_c\}} \\
		&\cdot f^\lambda_{\lambda_1,\nu,\lambda_2}(v) \cdot \frac{a_{2\nu_f}(v)}{a_{\nu_f}(v^2)}b_{\nu_c}(v) \cdot \fu_{\varrho \lambda_2} \oplus \fu_{\lambda_1^\diamond}*\K_{\lambda_2+\nu} \\
		=&\sum_{\lambda_1,\lambda_2 \in \fp, \nu \in \cf} v^{\Xi} f^\lambda_{\lambda_1,\nu,\lambda_2}(v) \cdot \frac{a_{2\nu_f}(v)}{a_{\nu_f}(v^2)}b_{\nu_c}(v) \cdot\K_{\lambda_2+\nu} \diamond \fU_{\varrho \lambda_2} \oplus \fU_{\lambda_1^\diamond}.
	\end{align*}
	where
	\begin{align*}
		\Xi=& \dim_\bfk\End(M_q(\lambda_1))+\dim_\bfk\End(M_q(\lambda_2))-\dim_\bfk\End(M_q(\lambda))\\
		&+\EuFo{\lambda_1,\nu+\lambda_2}_Q+\EuFo{\nu,\lambda_2}_Q+\frac{1}{2}\EuFo{\nu,\nu}_Q+\sum_{i \in \I_0} d^{\lambda_2}_i+\frac{1}{2}\sum_{i \in \I_0} d^\nu_i+2\{\nu_f,\nu_c\}.
	\end{align*}
	
	For type A, any $\beta \in \Phi^+$ could be represented by a subinterval $[s,t]$ of $[1,n]$, that is, $\beta=\alpha_s+\alpha_{s+1}+\cdots+\alpha_t$. We denote by $\fu_\beta$ (resp. $\fU_\beta$) by $\fu_{[s,t]}$ (resp. $\fU_{[s,t]}$). For convenience,  $\fu_{[s,t]}$ (and also $\fU_{[s,t]})$ represents the unit of $\widetilde{\ch}(Q,\btau)$ if $s>t$. Similar for $\K_{[s,t]}$.
	
	\begin{example}
		\label{ex: positivity of dcb of split A}
		Let $Q$ be the quiver of type $A_n$:
		\[\xymatrix{1\ar[r]&2\ar[r]&\cdots\ar[r]&n-1\ar[r]&n }\]
		and $\btau=\Id$. 
		For $\fU_{[s,t]}\in\widetilde{\ch}(Q,\btau)$, we have
		\begin{align*}
			\widetilde{\Omega}(\fU_{[s,t]})=&\sum_{k=-1}^{t-s} v^{\dim_\bfk\End(M_q([s,s+k]))+\dim_\bfk\End(M_q([s+k+1,t]))-1+\EuFo{[s,s+k],[s+k+1,t]}_Q+t-s-k}\\
			&\qquad \cdot f^{[s,t]}_{[s,s+k],[s+k+1,t]}(v) \cdot \K_{[s+k+1,t]} \diamond \fU_{[s+k+1,t]} \oplus \fU_{[s,s+k]^\diamond} \\
			=&\sum_{k=-1}^{t-s} v^{t-s-k} 
			\cdot \K_{[s+k+1,t]} \diamond \fU_{[s+k+1,t]} \oplus \fU_{[s,s+k]^\diamond},
		\end{align*}
		for any $1 \le s \le t \le n$. Note that $\fU_{[s+k+1,t]}$ and $\fU_{[s,s+k]^\diamond}$ commute for any $-1 \le k \le t-s$, and 
		$$\fU_{[s+k+1,t]}*\fU_{[s,s+k]^\diamond}=\fU_{[s+k+1,t]} \oplus \fU_{[s,s+k]^\diamond}.$$
		So $\fU_{[s+k+1,t]} \oplus \fU_{[s,s+k]^\diamond}$ is bar-invariant, and then
		$\mathfrak{L}_{[s+k+1,t],[s,s+k]^\diamond}=\fU_{[s+k+1,t]} \oplus \fU_{[s,s+k]^\diamond}$ by using Theorem \ref{iHA dCB theorem}. Therefore, the above formula is already a linear combination of dual canonical bases of $\widetilde{\ch}(Q^{\dbl},\swa)$, and all coefficients are in $\N[v^{\frac12},v^{-\frac12}]$.
	\end{example}

	\begin{example}
		\label{ex: positivity of dcb of quasi split A}
		We consider the quasi-split iquantum group of type $(A_{2n+1},\varrho)$ with $n\geq1$, $\varrho\neq\Id$. 
		Fix the iquiver given in \eqref{diag: A}. Let us compute $\widetilde{\Omega}(\fU_\beta)$ in $\widetilde{\ch}(Q,\btau)$, for any $\beta\in\Phi^+$. 
		
		By induction on $n$, one only needs to consider $\widetilde{\Omega}(\fU_{[1,m]})$ for $1 \le m \le 2n+1$ by symmetry. The computation is divided into the following two cases. 
		
		(1) $ \underline{1 \le m \le n+1}$. In this case, we have
		\begin{align*}
			\widetilde{\Omega}(\fU_{[1,m]})=&\sum_{k=0}^{m}v^{\dim_\bfk\End(M_q([1,k]))+\dim_\bfk\End(M_q([k+1,m]))-1+\EuFo{[1,k],[k+1,m]}_Q+(1-\delta_{k,m})\delta_{m,n+1}} \\
			& \cdot\K_{[k+1,m]} \diamond \fU_{[2n+2-m,2n+1-k]} \oplus \fU_{[1,k]^\diamond}\\
			=& \sum_{k=0}^{m} v^{\delta_{m,n+1}(1-\delta_{k,m})} \cdot \K_{[k+1,m]} \diamond \fU_{[2n+2-m,2n+1-k]} \oplus \fU_{[1,k]^\diamond}.
		\end{align*}
		Since $\fU_{[2n+2-m,2n+1-k]}$ and $\fU_{[1,k]^\diamond}$ commute and
		\[ \fU_{[2n+2-m,2n+1-k]} *\fU_{[1,k]^\diamond}= \fU_{[2n+2-m,2n+1-k]} \oplus \fU_{[1,k]^\diamond},\]
		one gets $\mathfrak{L}_{[2n+2-m,2n+1-k],1,k]^\diamond}=\fU_{[2n+2-m,2n+1-k]} \oplus \fU_{[1,k]^\diamond}$. Therefore, $\widetilde{\Omega}(\fU_{[1,m]})$ is a linear combination of dual canonical bases of $\widetilde{\ch}(Q^{\dbl},\swa)$ with all coefficients in $\N[v^{\frac{1}{2}},v^{-\frac{1}{2}}]$.

		(2) $\underline{n+2 \le m \le 2n+1}$. In this case, we have
		\begin{align*}
			&\widetilde{\Omega}(\fU_{[1,m]})=\fU_{[1,m]^\diamond}+\sum_{k=0}^n \sum_{l=0}^{m-n-1} v \cdot \K_{[n+1-k,n+1+l]} \diamond \fU_{[n+1-l,n+1+k]} \oplus \fU_{[1,n-k]^\diamond \oplus [n+2+l,m]^\diamond} \\
			+&\sum_{k=0}^{m-n-2}\sum_{l=1}^{m-n-1-k} (v^2-1) \cdot \K_{[n+1-k-l,n+1+k+l]} \diamond \fU_{[n+1-k,n+1+k]} \oplus \fU_{[1,n-k-l]^\diamond \oplus [n+2+k+l,m]^\diamond}.
		\end{align*}
		Here and below, $\fU_{[1,n-k-l]^\diamond \oplus [n+2+k+l,m]^\diamond}\in \widetilde{\ch}(Q^{\dbl},\swa)$ corresponds to the module $$[M_q(\alpha_1+\cdots+\alpha_{n-k-l})^\diamond\oplus M_q(\alpha_{n+2+k+l}+\cdots+\alpha_m)^\diamond]\in \widetilde{\ch}(\bfk Q^{\dbl},\swa).$$
		For $0 \le k \le m-n-2$ and $1 \le l \le m-n-1-k$, 
		similarly to Case (1), we can see \begin{align}
			\mathfrak{L}_{[n+1-k,n+1+k],[1,n-k-l]^\diamond \oplus [n+2+k+l,m]^\diamond}=\fU_{[n+1-k,n+1+k]} \oplus \fU_{[1,n-k-l]^\diamond \oplus [n+2+k+l,m]^\diamond}.
		\end{align}
		
		Now we focus on the terms $\fU_{[n+1-l,n+1+k]} \oplus \fU_{[1,n-k]^\diamond \oplus [n+2+l,m]^\diamond}$ for $0 \le k \le n$ and $0 \le l \le m-n-1$. Note that 
		\[ \fU_{[1,n-k]^\diamond \oplus [n+2+l,m]^\diamond}* \fU_{[n+1-l,n+1+k]}=\fU_{[n+1-l,n+1+k]} \oplus \fU_{[1,n-k]^\diamond \oplus [n+2+l,m]^\diamond}, \]
		since $n+1-l \le n+1 \le n+1+k$. We compute it in the following three subcases.
		
		(2a) \underline{$k=l$}. Since $\fU_{[n+1-l,n+1+k]}$ and $\fU_{[1,n-k]^\diamond \oplus [n+2+l,m]^\diamond}$ commute, we see 
		\[\mathfrak{L}_{[n+1-l,n+1+k],[1,n-k]^\diamond \oplus [n+2+l,m]^\diamond}=\fU_{[1,n-k]^\diamond \oplus [n+2+l,m]^\diamond}.\]
		
		(2b) \underline{$k<l$}. Then $n+1-l \le n-k$ and $n+1+k < n+2+l$. A simple computation shows 
		\begin{align*}
			&\fU_{[n+1-l,n+1+k]}* \fU_{[1,n-k]^\diamond \oplus [n+2+l,m]^\diamond}
			\\
			&=v^{-\frac{1}{2}(3-\delta_{k,n}-\delta_{l,m-n-1})} \fu_{[n+1-l,n+1+k]}* \fu_{[1,n-k]^\diamond \oplus [n+2+l,m]^\diamond}\\
			&=v^{-\frac{1}{2}(3-\delta_{k,n}-\delta_{l,m-n-1})}(\fu_{[n+1-l,n+1+k]}\oplus \fu_{[1,n-k]^\diamond \oplus [n+2+l,m]^\diamond}\\
			&\quad+v^{\delta_{l,n}}(v-v^{-1})\fu_{[n+1-k,n+1+k]} \oplus \fu_{[1,n-l]^\diamond \oplus [n+2+l,m]^\diamond} *\K_{[n+1-l,n-k]})  \\
			&=\fU_{[n+1-l,n+1+k]}\oplus \fU_{[1,n-k]^\diamond \oplus [n+2+l,m]^\diamond} \\
			&\quad+ (v-v^{-1}) \K_{[n+1-l,n-k]} \diamond \fU_{[n+1-k,n+1+k]} \oplus \fU_{[1,n-l]^\diamond \oplus [n+2+l,m]^\diamond}.
		\end{align*}
		Since $\mathfrak{L}_{[n+1-k,n+1+k],[1,n-l]^\diamond \oplus [n+2+l,m]^\diamond}=\fU_{[n+1-k,n+1+k]} \oplus \fU_{[1,n-l]^\diamond \oplus [n+2+l,m]^\diamond}$, we have
		\begin{align*}
			\mathfrak{L}_{[n+1-l,n+1+k],[1,n-k]^\diamond \oplus [n+2+l,m]^\diamond}=&\fU_{[n+1-l,n+1+k]}\oplus \fU_{[1,n-k]^\diamond \oplus [n+2+l,m]^\diamond}\\
			&-v^{-1} \cdot \K_{[n+1-l,n-k]} \diamond \mathfrak{L}_{[n+1-k,n+1+k],[1,n-l]^\diamond \oplus [n+2+l,m]^\diamond}.
		\end{align*}
		
		(2c)  \underline{$k>l$}. A similar argument shows that
		\begin{align*}
			\mathfrak{L}_{[n+1-l,n+1+k],[1,n-k]^\diamond \oplus [n+2+l,m]^\diamond}=&\fU_{[n+1-l,n+1+k]}\oplus \fU_{[1,n-k]^\diamond \oplus [n+2+l,m]^\diamond}\\
			&-v^{-1}\cdot \K_{[n+2+l,n+1+k]} \diamond \mathfrak{L}_{[n+1-l,n+1+l],[1,n-k]^\diamond \oplus [n+2+k,m]^\diamond}.
		\end{align*}
		
		Therefore, we have
		\begin{align*}
			&\widetilde{\Omega}(\fU_{[1,m]})\\
			=&\fU_{[1,m]^\diamond}+\sum_{k=0}^n \sum_{l=0}^{m-n-1} v \cdot \K_{[n+1-k,n+1+l]} \diamond \fU_{[n+1-l,n+1+k]} \oplus \fU_{[1,n-k]^\diamond \oplus [n+2+l,m]^\diamond} \\
			&+\sum_{k=0}^{m-n-2}\sum_{l=1}^{m-n-1-k} (v^2-1) \cdot \K_{[n+1-k-l,n+1+k+l]} \diamond \fU_{[n+1-k,n+1+k]} \oplus \fU_{[1,n-k-l]^\diamond \oplus [n+2+k+l,m]^\diamond} \\
			=& \fU_{[1,m]^\diamond}+\sum_{k=0}^n \sum_{l=0}^{m-n-1} v \cdot \K_{[n+1-k,n+1+l]} \diamond\mathfrak{L}_{[n+1-l,n+1+k],[1,n-k]^\diamond \oplus [n+2+l,m]^\diamond} \\
			&+ \sum_{0 \le k < l \le m-n-1} \K_{[n+1-l,n+1+l]} \diamond \mathfrak{L}_{[n+1-k,n+1+k],[1,n-l]^\diamond \oplus [n+2+l,m]^\diamond} \\
			&+\sum_{\substack{0 \le l < k \le n \\ l \le m-n-1}} \K_{[n+1-k,n+1+k]} \diamond \mathfrak{L}_{[n+1-l,n+1+l],[1,n-k]^\diamond \oplus [n+2+k,m]^\diamond} \\
			&+\sum_{0 \le k < l \le m-n-1} (v^2-1) \cdot \K_{[n+1-l,n+1+l]} \diamond \mathfrak{L}_{[n+1-k,n+1+k],[1,n-l]^\diamond \oplus [n+2+l,m]^\diamond} \\
			=&\fU_{[1,m]^\diamond}+\sum_{k=0}^n \sum_{l=0}^{m-n-1} v \cdot \K_{[n+1-k,n+1+l]} \diamond\mathfrak{L}_{[n+1-l,n+1+k],[1,n-k]^\diamond \oplus [n+2+l,m]^\diamond} \\
			&+ \sum_{0 \le k < l \le m-n-1} v^2 \cdot \K_{[n+1-l,n+1+l]} \diamond \mathfrak{L}_{[n+1-k,n+1+k],[1,n-l]^\diamond \oplus [n+2+l,m]^\diamond} \\
			&+\sum_{\substack{0 \le l < k \le n \\ l \le m-n-1}} \K_{[n+1-k,n+1+k]} \diamond \mathfrak{L}_{[n+1-l,n+1+l],[1,n-k]^\diamond \oplus [n+2+k,m]^\diamond}.
		\end{align*}
		It is clear that all coefficients are in $\N[v^{\frac{1}{2}},v^{-\frac{1}{2}}]$.
	\end{example}
	
	Finally, let us give two remarks.
	\begin{remark}
		From \cite{CZ26} and also the two examples, we find that it is difficult to prove Conjecture \ref{conj:positive} by direct computation, even for some special dual canonical basis elements. The potential strategy is to develop the geometric realization of quantum groups and iquantum groups in \cite{Qin,LW21b} further, in particular, try to interpret the comultiplication and embedding $\imath:\tUi\to \tU$ in the geometric setting. This strategy is natural and standard in geometric representation theory; cf. \cite{Lus90,Lus93,Na01,Qin,FL21,LW21b}.  
	\end{remark}
	
	\begin{remark}
		For the iquantum group $\tUi$ of type ${\rm AIII}_{2n}$ ($n\geq1$), its  Hall algebra realization has not been given; cf. Table \ref{tab:Satakediag}. However, 
		the dual canonical bases are defined for $\tUi$ and $\tU$ of arbitrary finite type in \cite{CLPRW}.  We further conjecture that the dual canonical basis of $\tUi$ of type ${\rm AIII}_{2n}$ has the properties outlined in Conjecture \ref{conj:positive}; moreover, the structure constants of its multiplication are also in $\N[v^{\frac12},v^{-\frac12}]$.
	\end{remark}	
	
	%
	
	%

	%


\end{document}